%% file: CarlementArXiv.tex
\documentclass[hidelinks,onefignum,onetabnum]{siamart220329}

\input{ex_shared}

\ifpdf
\hypersetup{
  pdftitle={Carleman Linearization of Nonlinear  Systems  and Its Finite-Section Approximations},
  pdfauthor={A. Amini, C. Zheng, Q. Sun and N. Motee}
}
\fi

\begin{document}

\maketitle

\begin{abstract}
The Carleman linearization is one of the mainstream approaches to lift a \linebreak[4] finite-dimensional nonlinear dynamical system
into an infinite-dimensional linear system with the promise of providing accurate approximations of the original nonlinear system over larger regions around the equilibrium for  longer time horizons with respect to the  conventional first-order linearization approach. Finite-section approximations of the lifted system has been widely used to study dynamical and control properties of the original nonlinear system. In this context, some of the outstanding problems are to determine under what conditions, as the finite-section order (i.e., truncation length) increases, the trajectory of the resulting approximate linear system from the finite-section scheme converges  to that of the original nonlinear system and whether the time interval over which the convergence happens can be quantified explicitly. In this paper, we provide  explicit error bounds for the finite-section approximation and prove that the convergence is indeed exponential with respect to the finite-section order. For a class of nonlinear systems, it is shown that  one can achieve exponential convergence over the entire time horizon up to infinity. Our results are practically plausible as our proposed error bound estimates can be used to compute proper truncation lengths for a given application, e.g., determining proper sampling period
for model predictive control and reachability analysis for safety verifications. We validate our theoretical findings through several illustrative simulations.
\end{abstract}



\section{Introduction}

For decades, the time-varying and nonlinear nature of most natural, physical, and engineered systems have imposed fundamental challenges on researchers to devise efficient and tractable algorithms to analyze and design such systems. In this work, we are interested in the class of time-varying nonlinear  systems whose dynamics are governed by
\begin{equation}\label{dynamicsystem}
	\dot{\bf x}(t)	
	= {\bf f}(t, {\bf x}(t))
\end{equation}
for all $t\ge t_0$ and  ${\bf x}(t_0)= {\bf x}_0\ne {\bf 0}$, with  the origin ${\bf 0}$ as their equilibrium, where  ${\bf x}\in {\mathbb R}^d$
is the state of the system and  
${\bf f}(t, {\bf x})$ is  an analytic function about ${\bf x}$ on a neighborhood of the equilibrium  ${\bf 0}$. A traditional method to study    system \eqref{dynamicsystem} is the well-known first-order linearization approach that   relies on obtaining a linear system  from the  first-order approximation of the function ${\bf f}(t, {\bf x})$. The first-order linearization approach is  valid only when the system is operating near its working point and over a short period of time. To study  various properties of  nonlinear  dynamical systems,  researchers have developed several frameworks over the past century \cite{  arnold1988dynamical,
	isidori1985nonlinear, khalil2015nonlinear, 
	kowalski1991nonlinear, wiggins2003introduction}.
One of  mainstream approaches  
is to lift the finite-dimensional nonlinear system \eqref{dynamicsystem} into an infinite-dimensional linear system.
Carleman linearization  and Koopman operator  are two of the most prominent methods that are closely connected in spirit
\cite{amini2019carleman, brockett2014early,   forets2017explicit,  forets2021reachability,  kowalski1991nonlinear, liu2021efficient,  minisini2007carleman,   steeb1980non}.
In this paper, we consider  Carleman linearization of the nonlinear dynamical system \eqref{dynamicsystem}, quantify several error bounds for its finite-section   approximations, and show that the resulting linear systems, for large enough truncation lengths,  provide precise approximation to the original nonlinear system on larger neighborhoods around the equilibrium, with respect to the first-order linearization, and for longer periods of time.

The control system community has experienced several success stories through methods that are developed based on Carleman linearization ideas \cite{amini2020approximate, amini2020quadratization,  forets2021reachability, hashemian2015fast, loparo1978estimating, pruekprasert2022moment, rauh2009carleman,krener1974linearization,krener1975bilinear, rotondo2022towards}.  For instance, the author of  \cite{rauh2009carleman} identified some connections between Carleman linearization and Lie series, and  then utilized it to design optimal control laws for infinite-dimensional systems.  Reference   \cite{loparo1978estimating} exploits Carleman approximation to obtain a relation between the lifted system and the domain of attraction of the original nonlinear system.  Recent work  \cite{hashemian2015fast}  employs Carleman linearization to implement model predictive control for nonlinear systems efficiently. In \cite{minisini2007carleman},  ideas from Carleman linearization are applied for state estimation and  design of feedback control laws. By exploiting the inherent structure of the lifted system, the authors of  \cite{amini2020approximate, amini2020quadratization} proposed a tractable method to quantilize and solve the Hamilton-Jacobi-Bellman equation through an exact iterative method.

The Carleman linearization  of the nonlinear  dynamical system \eqref{dynamicsystem}, which is expressed in  \eqref{carleman.eq3}, is an infinite-dimensional linear time-varying system, whose state matrix  is an upper-triangular block matrix. One should be meticulous in handling   the resulting infinite-dimensional linear system since the corresponding state matrix does {\it not} represent a bounded operator on the Hilbert space of all square-summable sequences. Moreover, the initial of the linear system has exponential decay when $\|{\bf x}_0\|_\infty<1$ and
exponential growth when
$\|{\bf x}_0\|_\infty>1$, where $\|{\bf x}_0\|_\infty$ is the maximal norm of the initial ${\bf x}_0$. These factors have prevented us from directly applying the existing theory  to analyze the resulting linear system from Carleman linearization  and then better understand the original nonlinear dynamical system.

A common remedy to deal with  the  Carleman linearization is finite-section approximations, which is given by \eqref{finitesection-1}, where one truncates the infinite-dimensional linear system   \cite{forets2017explicit, liu2021efficient,pruekprasert2022moment}.
A fundamental question is whether the first block of the solution
of the finite-section approximation converges to the solution of the original nonlinear dynamical system and what the convergence rate is.
In this work,  we provide a partial answer to the above questions when the coefficients in Maclaurin  expansion of the analytic function $f(t, {\bf x})$ enjoys certain uniform decay property, see Assumption \ref{assump-1}. In our main contribution, it is shown that
if the initial condition is in a vicinity of  the equilibrium, then
the first block of the solution  of the finite-section approximation will
exponentially converge  to the solution of the original nonlinear  system
as the order of the truncation in the finite-section approximation increases, see Theorems  \ref{maintheorem}, \ref{maintheorem2} and \ref{maintheorem3}.
We highlight that  the authors of   \cite{forets2017explicit, liu2021efficient} have established similar convergence result when the  function $f(t,{\bf x})$  in    \eqref{dynamicsystem} is a
 polynomial.

The paper is organized as follows. In Section \ref{introduction.section},  we consider
Carleman linearization of the nonlinear dynamical system  \eqref{dynamicsystem} and its finite section approximation, see \eqref{carleman.eq3} and \eqref{finitesection-1}.
In Section \ref{convergence.section},
we establish the exponential convergence of the  finite section scheme  over some time interval, see Theorem  \ref{maintheorem},
Corollary \ref{maincorollary.polynomial},
and  Theorems \ref{maintheorem2} and  \ref{maintheorem3}.
In Section \ref{numericalsimulation.section},  Carleman linearization of
several benchmark systems are discussed to  validate and illustrate the theoretical findings in Theorems  \ref{maintheorem} and \ref{maintheorem2}. 
The technical proofs of all conclusions 
are provided in Section \ref{proof.section}.

Some preliminary versions of this work were announced in \cite{amini2019carleman, amini2021error}. The authors assert that the content of this manuscript  significantly differs from its conference versions as this work contains several new  and improved results as well as several new case studies with respect to its conference versions.

\smallskip

\section{Carleman Linearization and Its Finite-Section Approximations}\label{introduction.section}

We review some notions related to the Carleman linearization
and its finite section scheme
\cite{amini2021error, amini2022error, forets2017explicit, forets2021reachability, liu2021efficient,
 minisini2007carleman,  pruekprasert2022moment,rauh2009carleman}. For a  given vector  ${\bf x}=[x_1, \ldots, x_d]^T\in {\mathbb R}^d$,  let us denote monomial
${\bf x}_{\pmb \alpha}= x_1^{\alpha_1} \cdots x_d^{\alpha_d}$ for some  ${\pmb \alpha}=[\alpha_1, \ldots, \alpha_d]^T\in {\mathbb Z}_+^d$, where ${\mathbb Z}_+^d$ is the set of all $d$-dimensional non-negative integer vectors.
Suppose that  the Maclaurin series of the vector-valued analytic function ${\bf f}(t, {\bf x})$ in \eqref{dynamicsystem}  can be expressed as
\begin{equation}
	\label{marcluurin.def}
	{\bf f}(t, {\bf x})=\sum_{{\pmb \alpha}\in { \mathbb Z}_+^d\backslash \{{\bf 0}\}} {\bf f}_{\pmb\alpha}(t) {\bf x}_{\pmb\alpha}
\end{equation}
for all $t\ge t_0$. For a given ${\pmb \alpha}=[\alpha_1, \ldots, \alpha_d]^T\in {\mathbb Z}_+^d$, let us define $|\pmb\alpha|=\alpha_1+\cdots+\alpha_d$  and set ${ \mathbb Z}_k^d=\{ {\pmb \alpha}\in {\mathbb Z}_+^d~ |~  |\pmb \alpha|=k \}$. In this work, we  assume that the coefficients of the Maclaurin series \eqref{marcluurin.def} have the following
uniform exponential decay property; we refer to    \eqref{maccoefficient.remark.eq1} for another conventional uniform exponential decay assumption.

\begin{assum}\label{assump-1}  There exist positive constants $D_0$ and $R>0$  such that the coefficients ${\bf f}_{\pmb \alpha}(t)=[f_{1, \pmb \alpha}(t), \ldots, f_{d, \pmb \alpha}(t)]^T$ in the
	Maclaurin expansion \eqref{marcluurin.def} satisfy
	\begin{equation}  \label{fassump.eq1}
		\sup_{t\ge t_0} \sum_{j=1}^d \sum_{\pmb \alpha\in { \mathbb Z}_k^d} |f_{j, \pmb\alpha}(t)|\le  D_0 R^{-k}
	\end{equation}
	for $k\ge 1$.
\end{assum}

When   ${\bf f}(t,{\bf x})$ in  \eqref{dynamicsystem} is a  polynomial of degree $L\ge 1$, i.e.,
\begin{equation}\label{maincorollary.eq1}
{\bf f}(t,{\bf x})=	{\bf p}_L(t, {\bf x})=
	\sum_{1\le |\pmb \alpha|\le L} \hspace{0.05cm} [p_{1, \pmb \alpha}(t), \ldots,p_{d, \pmb \alpha}(t)]^T  \hspace{0.02cm}  {\bf x}^{\pmb \alpha},
\end{equation}
Assumption \ref{assump-1} will be  
satisfied as one can  verify that for every convergence radius $R>0$,
the uniform exponential decay property \eqref{fassump.eq1}  holds
with $D_0$ replaced by
\begin{equation}
	\label{maincorollary.eq2}
	D_0({\bf p}_L, R)= \sup_{1\le k\le L}  R^{k} \sup_{t\ge t_0} \sum_{j=1}^d \sum_{\pmb \alpha\in {\mathbb Z}_k^d} |p_{j, \pmb \alpha}(t)|. \end{equation}

\vspace{0.1cm}

Let us denote the standard Euclidean basis for $\R^d$ by ${\bf e}_1, \ldots, {\bf e}_d$  and set
\begin{equation}\label{fj.negative}
	f_{j, \pmb \alpha}(t)=0
\end{equation}
for all $\pmb \alpha\not\in {{\mathbb Z}}_+^d\backslash \{{\bf 0}\}$ and 	$j=1, \ldots, d$. The  Carleman linearization of the nonlinear  dynamical system \eqref{dynamicsystem} starts from its reformulation
\begin{equation}\label{dynamicsystem.reformulation}
	\dot{x}_j (t) 
	=\sum_{{\pmb \alpha}\in  { \mathbb Z}_+^d\backslash \{{\bf 0}\}} f_{j, {\pmb \alpha} }(t) \hspace{0.05cm}{\bf x}_{\pmb \alpha}(t)
\end{equation}
for $j=1,\ldots,d$.
For every $\pmb \alpha=[\alpha_1, \ldots, \alpha_d]\in { \mathbb Z}_+^d\backslash \{{\bf 0}\}$,
the derivative of  monomial $\xa$
can be calculated as
\begin{eqnarray}\label{carleman.eq0}
	\dot{\x}_{\a}(t)  
	& = &
	\sum_{j=1}^d  ~\alpha_j {\bf x}_{\pmb \alpha-{\bf e}_j} \sum_{\pmb \gamma \in { \mathbb Z}_+^d \backslash \{{\bf 0}\}} f_{j, \pmb \gamma}(t) \hspace{0.05cm} {\bf x}_{\pmb \gamma}(t) \\
	& = & \sum_{\pmb \beta\in { \mathbb Z}_+^d\backslash \{{\bf 0}\}} \left(\sum_{j=1}^d \alpha_j f_{j, \pmb \beta-\pmb \alpha+{\bf e}_j}(t)\right)  {\bf x}_{\pmb \beta}(t), \nonumber
\end{eqnarray}
with  initial condition
${\bf x}_{\pmb \alpha}(t_0)= ({\bf x}_0)_{\pmb \alpha}$. 
For every  $k\ge 1$, we define a new state variable as  ${\bf z}_k=[{\bf x}_{\pmb \alpha}]_{{\pmb \alpha}\in {\mathbb Z}_k^d}$, which contains all the monomials of order $k$.
Regrouping monomials   in \eqref{carleman.eq0} all  together yields the following infinite-dimensional linear  system
\begin{equation}\label{carleman.eq1}
	\dot{\bf z}_k(t) = \sum_{l=k}^\infty {\bf A}_{k, l}(t) \hspace{0.05cm}{\bf z}_l(t)  \ \ {\rm and} \ \  {\bf z}_k(t_0)=\big[{\bf x}^0_{\pmb \alpha}\big]_{{\pmb \alpha}\in {\mathbb Z}_k^d}
\end{equation}
for all  $t\ge t_0$  and $k\ge 1$, 
where
\begin{equation}  \label{carleman.eq2}
	{\bf A}_{k, l}(t)= \left[\sum_{j=1}^d \alpha_j f_{j, \pmb\beta-\pmb\alpha+{\bf e}_j}(t)\right]_{{\pmb \alpha}\in {\mathbb Z}_k^d,~ {\pmb \beta}\in {\mathbb Z}_l^d}
\end{equation}
are matrices of size $ { {k+d-1}\choose d-1}\times { {l+d-1}\choose d-1}$.
By defining the infinite-dimensional state vector
\begin{equation}\label{z.def} {\bf z}=[{\bf z}_1, {\bf z}_2, \ldots,  {\bf z}_N, \ldots]^T,
\end{equation}  the set of linear systems \eqref{carleman.eq1}  can be rewritten in the following infinite-dimensional matrix form
\begin{equation}  \label{carleman.eq3}
	\dot{\bf z}(t)= {\bf A}(t) {\bf z}(t)
\end{equation}
for all $t\ge t_0$ with the initial condition ${\bf z}(t_0)=[{\bf z}_k(t_0)]_{k\ge 1}$, where  
\begin{equation}\label{At.def}
	{\bf A}(t)= \left[\begin{array}{ccccc}
		{\bf A}_{1, 1}(t) & {\bf A}_{1, 2}(t) & \cdots  & {\bf A}_{1, N}(t) & \cdots\\
		& {\bf A}_{2, 2}(t)& \cdots & {\bf A}_{2, N}(t) & \cdots\\
		& &\ddots & \vdots & \ddots\\
		& & & {\bf A}_{N, N}(t) & \cdots\\
		& & & & \ddots
	\end{array}\right].\end{equation}
The resulting linear system  \eqref{carleman.eq3}  is referred  to   as {\it Carleman linearization} of the  nonlinear dynamical system \eqref{dynamicsystem}.   While
the state-space of the original nonlinear system \eqref{dynamicsystem} is the finite-dimensional Euclidean space $\R^d$, its Carleman linearization \eqref{carleman.eq3}  is an infinite-dimensional  linear time-varying system
whose state matrix ${\bf A}(t)$ is an   { upper-triangular block matrix}.
According to the bound estimates for the block matrices ${\bf A}_{k,l}(t)$ for all  $1\le k\le l$ in Lemma \ref{maintheorem.lem1},
the  state matrix ${\bf A}(t)$ in \eqref{carleman.eq3} is not  a bounded operator on $\ell^2({\mathbb Z}_+^d\backslash \{{\bf 0}\})$, the Hilbert space of all square-summable sequences on
${\mathbb Z}_+^d\backslash \{{\bf 0}\}$.
Moreover, it is observed that
its initial  ${\bf z}(t_0)$  has exponential decay when $\|{\bf x}_0\|_\infty<1$ and exponential growth when
$\|{\bf x}_0\|_\infty>1$, in which the maximal norm of
a vector ${\bf x}=[x_1, \ldots, x_d]^T\in {\mathbb R}^d$ is represented 	by $\|{\bf x}\|_\infty=\max_{1\le j\le d} |x_j|$.
The above two observations  have been the main preventive factors to apply existing theory on Hilbert space directly to analyze the
Carleman linearization of a nonlinear  system.

\smallskip

A conventional approach to solve the infinite-dimensional  linear system \eqref{carleman.eq3} is to consider its finite-section approximation of order $N$, which is given by
\begin{equation}
	\label{finitesection-1}
	\left [\begin{array}{c}
		\dot{\bf y}_{1, N}(t)\\
		\dot{\bf y}_{2, N}(t)\\
		\vdots\\
		\dot{\bf y}_{N, N}(t)
	\end{array}\right]
	= \left[\begin{array}{ccccc}
		{\bf A}_{1, 1}(t) & {\bf A}_{1, 2}(t) & \cdots  & {\bf A}_{1, N}(t) \\
		& {\bf A}_{2, 2}(t)& \cdots & {\bf A}_{2, N}(t) \\
		& & \ddots & \vdots \\
		& & & {\bf A}_{N, N}(t)
	\end{array}\right]
	\left[\begin{array}{c}
		{\bf y}_{1, N}(t)\\
		{\bf y}_{2, N}(t)\\
		\vdots\\
		{\bf y}_{N, N}(t)
	\end{array}\right]
\end{equation}
with initial
${\bf y}_{k,N}(t_0)={\bf z}_k(t_0)$ for $k =1, \ldots, N$ \cite{ forets2017explicit,  liu2021efficient, mezic2022numerical, minisini2007carleman}.
The above finite-section scheme is of dimension
$ { {N+d}\choose d}-1$ and can be solved
 by first solving
$$ \dot{\bf y}_{N, N}(t)= {\bf A}_{N, N}(t) {\bf y}_{N, N}(t) \ \ {\rm with} \ \
y_{N,N}(t_0)={\bf z}_N(t_0) $$ 
and then solving
$$ \dot{\bf y}_{k, N}(t)= {\bf A}_{k, k}(t) {\bf y}_{k, N}(t)+ \sum_{l=k+1}^N {\bf A}_{k, l}(t) {\bf y}_{l, N}(t)
\ \ {\rm with} \ \
y_{k,N}(t_0)= {\bf z}_k(t_0) $$
for $k=N-1, \ldots, 1$,  recursively, where the size of each subsystem is ${{k+d-1} \choose {d-1}}$.
For $N=1$, the finite-section approximation
\eqref{finitesection-1} becomes the well-known first-order linearization of the nonlinear system  \eqref{dynamicsystem} that has been widely used in practice. The first block of the state vector in the linear system \eqref{carleman.eq3} corresponds to the solution of the original nonlinear system. When considering the finite-section approximation \eqref{finitesection-1}, it is desirous to scrutinize whether the first block of its solution converges to the solution of the original nonlinear system  \eqref{dynamicsystem}
and what the convergence rate is. In the next section,  we provide   partial answers to the above query when the function $f(t, {\bf x})$ in  \eqref{dynamicsystem} satisfies  Assumption \ref{assump-1};
c.f. \cite{forets2017explicit,  liu2021efficient} where
$f(t, {\bf x})$ is a polynomial.

\smallskip

\section{Convergence of Finite-Sectioning of the Carleman Linearization}
\label{convergence.section}

In this part, we  study convergence properties of    the finite-section approximation \eqref{finitesection-1}
of the Carleman linearization \eqref{carleman.eq3}.
Theorem 3.1 shows that if the initial condition ${\bf x}_0$ satisfies
\begin{equation}\label{maintheorem.eq2}
	0<\|{\bf x}_0\|_\infty <R/e,
\end{equation}
then the first block  ${\bf y}_{1, N}(t)$ for  $N\ge 1$, of the solution of the finite-section approximation \eqref{finitesection-1}
  converges  exponentially to  the true solution   ${\bf x} (t)$ of  the nonlinear dynamical  system \eqref{dynamicsystem}
over time interval $[t_0, t_0+T^*)$, where
$D_0$ and $R$ are the constants in Assumption \ref{assump-1} and
\begin{equation}\label{Tstar.def}
	T^*= \frac{(e-1)R}{(2e-1)D_0} \ln\left(\frac{R}{e\|{\bf x}_0\|_\infty}\right).
\end{equation}
In \cite[Theorem 4.2]{forets2017explicit}, the authors consider a special case of this problem when $f(t,{\bf x})$ is a time-independent polynomial. This implies that finite-section approximation of the Carleman linearization provides a reasonable approximation of the original nonlinear system over a quantifiable time interval that  depends linearly on the logarithmic scale of the
distance of the initial from the equilibrium.

Furthermore, we  are interested in those  nonlinear   systems \eqref{dynamicsystem} with the following additional property on their Jacobian $\nabla {\bf f}(t, {\bf 0})
=[{\bf f}_{{\bf e}_1}(t), \ldots, {\bf f}_{{\bf e}_d}(t)]^T$  
at the equilibrium  ${\bf 0}$.

\begin{assum}\label{assump-2}
	The Jacobian
	$\nabla {\bf f}(t, {\bf 0})$ is a  { time-independent} diagonal matrix with  negative diagonal entries
	$\lambda_{\pmb \alpha}$ for every $|\pmb \alpha|=1$ that satisfy
	\begin{equation} \label{fassump.eq2}
		\lambda_{\pmb \alpha} \le -\mu_0  ~~{\rm for \ all} ~~ |\pmb \alpha|=1
	\end{equation}
for some $\mu_0>0$.
\end{assum}

 In Theorem \ref{maintheorem2}, we show that if the initial condition ${\bf x}_0$ satisfies
\begin{equation}\label{maintheorem2.eq1}
	0<\|{\bf x}_0\|_2 < \frac{ R^2 \mu_0}{D_0+R \mu_0},
\end{equation}
where $\|{\bf x}_0\|_2$ is the Euclidean norm of $\x_0$,  then the first block  ${\bf y}_{1, N}$  of the solution of the finite-section approximation \eqref{finitesection-1}
of the Carleman linearization \eqref{carleman.eq3} will converge  exponentially to the true solution  of  the nonlinear dynamical  system \eqref{dynamicsystem} over the
 {\sl entire}  time interval $[t_0, \infty)$,
where $D_0, R, \mu_0$ are the constants in Assumptions \ref{assump-1} and \ref{assump-2}; cf.
\cite[Corollary 1 in the Supplementary Information]{liu2021efficient}
for the case that $f(t,{\bf x})$ is a  polynomial.
It is observed from  Corollary \ref{maintheorem3.cor1} that
the nonlinear dynamical system \eqref{dynamicsystem} with
the  analytic function ${\bf f}(t,{\bf x})$   satisfying
Assumptions \ref{assump-1} and \ref{assump-2} is stable when the  initial ${\bf x}(t_0)={\bf x}_0$ satisfies
\eqref{maintheorem2.eq1}. Despite the traditional first-order linearization approach that results in approximations that are useful  only over short time intervals, we
 demonstrate that
 the Carleman linearization can be employed over the entire time horizon when the origin is an asymptotically stable equilibrium of the nonlinear dynamical system.

\smallskip

The next theorem is the first main contribution of this work, where its proof is given in Section \ref{maintheorem.pfsection} and a couple of  supporting numerical examples are demonstrated in Section \ref{numericalsimulation.section}.

\begin{theorem}\label{maintheorem}
Suppose that ${\bf x}(t)$
 is the  solution of  the nonlinear dynamical  system \eqref{dynamicsystem}, the analytic function ${\bf f}(t,{\bf x})$ in  \eqref{dynamicsystem}  satisfies Assumption \ref{assump-1}, and ${\bf y}_{1, N}(t)$ for every $N\ge 1$  is the first block  of the solution of the finite-section approximation \eqref{finitesection-1}. If  the initial  ${\bf x}(t_0)={\bf x}_0$ satisfies \eqref{maintheorem.eq2},
	then     
	\begin{equation}\label{maintheorem.eq4}
		\|{\bf y}_{1, N}(t)-{\bf x} (t)\|_\infty\le  \frac{RM_0}{\sqrt{2\pi}(R-M_0)}
		N^{-3/2} e^{D_0(t-t_0)N/R} \left(\frac{\|{\bf x}_0\|_\infty e}{R}\right)^{ (e-1)N/(2e-1)}
	\end{equation}
	holds for all $t_0\le t\le t_0+T^*$ and $N\ge 1$,
	where  $T^*$ is defined by \eqref{Tstar.def} and
\begin{equation}\label{maintheorem.M0def}
		M_0=\|{\bf x}_0\|_\infty^{(e-1)/(2e-1)} (R/e)^{e/(2e-1)}< R/e.
	\end{equation}
\end{theorem}

An important step in the proof of Theorem \ref{maintheorem} is to establish the local bound estimate
\begin{equation} \label{maintheorem.lem2.eq1}
	\|{\bf x}(t)\|_\infty \le M_0  \ \ {\rm for \ all}\ \  t_0\le t\le t_0+T^*.
\end{equation}
We refer to Lemma \ref{maintheorem.lem2} for more details. Whenever a local bound for the solution ${\bf x}(t)$ for all $t\ge t_0$
is provided  a priori, i.e.,  there exist an upper bound $M<R/e$ and a time range $T_1>0$ such that
\begin{equation}\label{boundedness.assumption}
	\|{\bf x}(t)\|_\infty\le M\ \ {\rm for \ all}\ \ t_0\le t\le t_0+T_1,
\end{equation}
one can apply  a similar argument 
to establish the following result.

\begin{corollary}\label{main.cor}
Suppose that  the  solution  ${\bf x}(t)$   of  the nonlinear dynamical system \eqref{dynamicsystem}  satisfies
\eqref{boundedness.assumption}, and the analytic function ${\bf f}(t,{\bf x})$ in  \eqref{dynamicsystem}  satisfies Assumption \ref{assump-1}.
	Then  
	\begin{equation}\label{maincor.eq4}
		\|{\bf y}_{1, N}(t)-{\bf x} (t)\|_\infty\le \frac{RM}{\sqrt{2\pi}(R-M)}
		N^{-3/2}   \left(\frac{M\exp (D_0 (t-t_0)/R+1)}{R}\right)^N	
	\end{equation}
	holds for all $t_0\le t \le  t_0+ \tilde T^*$ and $N\ge 1$, where $M, D_0, R$  are constants in \eqref{boundedness.assumption} and Assumption \ref{assump-1} and
	\begin{equation}\label{Tstarcor.def}
		\tilde T^*=\min\left\{T_1,  \frac{R}{D_0} \ln\frac{R}{ Me}\right\}.
	\end{equation}
\end{corollary}

As a consequence of Theorem \ref{maintheorem}, one has the following result for systems \eqref{dynamicsystem} whose
right-hand side is a polynomial  described by
\eqref{maincorollary.eq1}; cf. \cite[Theorem 4.2]{forets2017explicit}.

\begin{corollary}
	\label{maincorollary.polynomial}
	Let us consider the nonlinear  system
	\eqref{dynamicsystem}   with a nonzero initial ${\bf x}(t_0)={\bf x}_0\ne 0$, where its  ${\bf f}(t, {\bf x})$ is a polynomial as in  \eqref{maincorollary.eq1} of order $L\ge 2$, and set
	\begin{equation} \label{maincorollary.eq3}
		T^*({\bf p}_L, {\bf x}_0)=\sup_{R> e\|{\bf x}_0\|_\infty} \frac{(e-1)R}{(2e-1)D_0({\bf p}_L, R)} \ln \left(\frac{R}{e\|{\bf x}_0\|_\infty}\right),\end{equation}
	where $	D_0({\bf p}_L, R)$ is given by \eqref{maincorollary.eq2}.
	Then, the first block   ${\bf y}_{1, N}(t)$
of the solution of the finite-section approximation \eqref{finitesection-1}
	converges exponentially
	to the  solution
	${\bf x}(t)$ of  the  original nonlinear  system
	\eqref{dynamicsystem}  for all $t_0<t < t_0+T^*({\bf p}_L, {\bf x}_0)$.
\end{corollary}

In order to calculate the maximal achievable time range from \eqref{maincorollary.eq3}, we define
$$a_k=\sup_{t\ge t_0} \sum_{j=1}^d \sum_{\pmb \alpha\in {\mathbb Z}_k^d} |p_{j, \pmb \alpha}(t)| $$
for all $k=1, \ldots, L$.  For all  $R\ge \max_{1\le k\le L-1} (a_k/a_L)^{1/(L-k)}$, one  may verify that
$D_0({\bf p}_L, R)  =   
a_L R^L$.
Hence,  whenever
$\|{\bf x}_0\|_\infty\ge  e^{-1} \max_{1\le k\le L-1} (a_k/a_L)^{1/(L-k)}$,
we should select $R\ge \max_{1\le k\le L-1} (a_k/a_L)^{1/(L-k)}$, which as a result the  maximal achievable time range according to \eqref{maincorollary.eq3} can be found as
\begin{eqnarray}
	T^*({\bf p}_L, {\bf x}_0) 
	& = & \frac{e-1}{(2e-1)a_L} \sup_{R> e \|{\bf x}_0\|_\infty } R^{1-L} \ln \left(\frac{R}{e \|{\bf x}_0\|_\infty} \right) \label{maincorollary.eq4} \\ 
	& = & \frac{e-1}{(2e-1)(L-1)e^L  a_L } \|{\bf x}_0\|_\infty^{1-L}. \nonumber
\end{eqnarray}
The last equality in \eqref{maincorollary.eq4} follows as it can be verified that $R^{1-L} (\ln  R- \ln  \|{\bf x}_0\|_\infty-1)$ takes its maximal values at $R= e^{L/(L-1)} \|{\bf x}_0\|_\infty$.
We remark that whenever  the polynomial
${\bf p}_{L}(t, {\bf x})$  in \eqref{maincorollary.eq1} for  $L\ge 2$
is time-independent,  the authors of \cite{forets2017explicit} provide an estimate for the maximal time range $T^*$
similar to the one in \eqref{maincorollary.eq4}; we refer to \cite[Theorem 4.3]{forets2017explicit}.

\smallskip

By Theorem \ref{maintheorem},
the first block  of the solution of the finite section scheme \eqref{finitesection-1}
enjoys exponential convergence to the true solution over a quantifable period of time.
In \cite{liu2021efficient},
  the authors consider
the  nonlinear  dynamical system
\eqref{dynamicsystem} with   ${\bf f}(t, {\bf x})$  being a  
polynomial
${\bf p}_{L}(t, {\bf x})$ as in \eqref{maincorollary.eq1} with $L=2$.
Under the assumptions that  the gradient
 $\nabla {\bf f}(t, {\bf 0})$  is time-independent, diagonalizable, and has eigenvalues that have negative real parts, c.f. Assumption
\ref{assump-2},  they show that
the first block    of the solution of the finite-section approximation \eqref{finitesection-1}
converges exponentially to the solution ${\bf x}(t)$ of   the nonlinear dynamical  system \eqref{dynamicsystem}
on the whole time range $t\ge t_0$
when the initial ${\bf x}_0$
is not too far away from the origin.
In the next theorem, which is our second main contribution of this work,   we consider exponential convergence of the first block  of the solution of  the finite-section approximation \eqref{finitesection-1}
over the entire time horizon from $t_0$ to $\infty$, where
${\bf f}(t, {\bf x})$ in  \eqref{dynamicsystem}
satisfies Assumptions \ref{assump-1} and \ref{assump-2}.

\begin{theorem} \label{maintheorem2}
	Suppose that
	${\bf x}(t)$ is the continuous solution of  nonlinear system \eqref{dynamicsystem}
	and  the  analytic function ${\bf f}(t,{\bf x})$ in   \eqref{dynamicsystem}  satisfy
	Assumptions \ref{assump-1} and \ref{assump-2}.
	If  the initial  ${\bf x}(t_0)={\bf x}_0$ satisfies  \eqref{maintheorem2.eq1},
	then
	\begin{equation} \label{maintheorem2.eq2}
		\|{\bf y}_{1, N}(t)-{\bf x}(t)\|_\infty \le
		\|{\bf x}_0\|_2 \left(\frac{ (D_0+R\mu_0)\|{\bf x}_0\|_2 }{ R^2 \mu_0}\right)^{N}
	\end{equation}
hold for all $N\ge 1$ and $t\ge t_0$, where $D_0, R, \mu_0$ are the constants in Assumptions \ref{assump-1} and \ref{assump-2}.
\end{theorem}

The proof of a stronger  version of Theorem \ref{maintheorem2}  is given in
Section \ref{maintheorem2.pfsection} and  some numerical demonstrations of the result of Theorem \ref{maintheorem2} is presented in Section \ref{numericalsimulation.section}.


In the conclusions of Theorems \ref{maintheorem} and \ref{maintheorem2}, it is assumed that the nonlinear dynamical system \eqref{dynamicsystem} enjoys the equilibrium assumption
${\bf f}(t, {\bf 0})= {\bf 0}$. In the next step, we generalize our results  in Theorem \ref{maintheorem2} to systems
\begin{equation} \label {nohomogenous.def}
\dot{\bf x}(t)={\bf f}(t, {\bf x})
\end{equation}
with  ${\bf x}(t_0)={\bf x}_0$,
 whose behavior at the origin, which is characterized by perturbation term ${\bf f}(t, {\bf 0})=[f_{1, {\bf 0}}(t),\cdots, f_{d, {\bf 0}}(t)]^T$,  satisfies
 \begin{equation}\label{assumption-3}
\sup_{t\ge t_0} \sum_{j=1}^d |f_{j, {\bf 0}}(t)|\le \nu_0
\end{equation}
for some small $\nu_0\ge 0$. For the analytic function ${\bf f}(t, {\bf x})$  satisfying \eqref{assumption-3} and
  Assumption \ref{assump-1},  we express its
 Maclaurin series  by
\begin{equation}
	\label{marcluurinnohomogeneous.def}
	{\bf f}(t, {\bf x})=\sum_{{\pmb \alpha}\in { \mathbb Z}_+^d} {\bf f}_{\pmb\alpha}(t) {\bf x}_{\pmb\alpha},
\end{equation}
where ${\bf f}_{\pmb\alpha}(t)=[f_{1, {\pmb \alpha}} (t), \ldots, f_{d, {\pmb \alpha}}(t)]^T$ are the coefficient of the series; c.f. \eqref{marcluurin.def} where $f_{\bf 0}(t)= f(t, {\bf 0})={\bf 0}$.
By applying the  procedure in Section \ref{introduction.section},
one can obtain the  Carleman linearization of the nonlinear dynamical system \eqref{dynamicsystem} as
\begin{equation}  \label{carlemannonhomo.eq3}
	\dot{\bf z}(t)= {\bf A}(t) {\bf z}(t)
\end{equation}
for all $t\ge t_0$ with the initial ${\bf z}(t_0)=[{\bf z}_k(t_0)]_{k\ge 1}$, where
${\bf z}(t)$ is
the infinite-dimensional state vector in \eqref{z.def}, block matrices ${\bf A}_{k, l}(t)$ for $k, l\ge 1$, are given in
\eqref{carleman.eq2},
and
\begin{equation*}
{\bf A}(t)
= \left[\begin{array}{ccccccc}
		{\bf A}_{1, 1}(t) & {\bf A}_{1, 2}(t) & {\bf A}_{1,3}(t) & \cdots  &  {\bf A}_{1, N-1}(t)  & {\bf A}_{1, N}(t) & \cdots\\
			  {\bf A}_{2, 1}(t)
		& {\bf A}_{2, 2}(t)& {\bf A}_{2, 3}(t) & \cdots & {\bf A}_{2, N-1}(t) & {\bf A}_{2, N}(t) & \cdots\\
		& {\bf A}_{3, 2}(t)& {\bf A}_{3, 3}(t) & \cdots &  {\bf A}_{3, N-1}(t) & {\bf A}_{2, N}(t) & \cdots\\
		& &\ddots & \ddots & \vdots &  \vdots & \ddots\\
	& 	& &\ddots   & {\bf A}_{N-1, N-1}(t) & {\bf A}_{N-1, N}(t)  & \ddots\\
		& & & &  {\bf A}_{N, N-1}(t) &  {\bf A}_{N, N}(t) & \ddots\\
		& & & & &  \ddots & \ddots
	\end{array}\right]. \end{equation*}
We emphasize the resulting state matrix in \eqref{carlemannonhomo.eq3} is no longer upper triangular.
 Similarly, one may verify that the finite-section approximation of the Carleman linearization
\eqref{carlemannonhomo.eq3} can be described as
 \begin{equation} \label{Athomogenous.def}
\hskip-0.02in	\left [\begin{array}{c}
		\dot{\bf y}_{1, N}\\
		\dot{\bf y}_{2, N}\\
		\vdots\\
\dot{\bf y}_{N-1, N}\\
		\dot{\bf y}_{N, N}
	\end{array}\right]
\hskip-0.02in = \hskip-0.02in \left[\begin{array}{cccccc}
		{\bf A}_{1, 1}(t) & {\bf A}_{1, 2}(t)  & \cdots  &  {\bf A}_{1, N-1}(t)  & {\bf A}_{1, N}(t) \\
			  {\bf A}_{2, 1}(t)
		& {\bf A}_{2, 2}(t)  & \cdots & {\bf A}_{2, N-1}(t) & {\bf A}_{2, N}(t) \\
		& \ddots 
& \ddots & \vdots & \vdots   \\
	 	& &\ddots   & {\bf A}_{N-1, N-1}(t) & {\bf A}_{N-1, N}(t)  \\
		& & &  {\bf A}_{N, N-1}(t) &  {\bf A}_{N, N}(t) \\
			\end{array}\right]
	\left[\begin{array}{c}
		{\bf y}_{1, N}\\
		{\bf y}_{2, N}\\
		\vdots\\
		{\bf y}_{N-1, N}\\
{\bf y}_{N, N}
	\end{array}\right] \end{equation}
where ${\bf y}_{k, N}:= {\bf y}_{k, N}(t)$ satisfies the initial condition ${\bf y}_{k, N}(t_0)={\bf z}_k(t_0)$.

Before stating a stronger  version of Theorem \ref{maintheorem2}, let us define parameters
 \begin{equation} \label{maintheorem3.eq1+}
    \eta_0:= \frac{\nu_0}{D_0} \ \ {\rm and} \ \   \eta_1:= \frac{R\mu_0}{D_0}\in (0, 1],
\end{equation}
where $D_0, R, \mu_0$ are constants in Assumptions \ref{assump-1} and \ref{assump-2}.

\begin{theorem}\label{maintheorem3}
Suppose that ${\bf x}(t)$ is the solution of the  nonlinear system \eqref{nohomogenous.def} whose
coefficients in the Maclaurin series \eqref{marcluurinnohomogeneous.def} satisfy   Assumptions \ref{assump-1} and \ref{assump-2} and
  \eqref{assumption-3}  for some
 $\nu_0\ge 0$. If
  \begin{equation}
 \label{maintheorem3.eq1}
 \eta_0\le  2+\eta_1-2\sqrt{1+\eta_1}
  \end{equation}
and  initial  ${\bf x}(t_0)={\bf x}_0$ satisfies
\begin{equation}   \label{maintheorem3.eq2}
\|{\bf x}_0\|_2< R\epsilon_0,
\end{equation}
then
\begin{equation} \label{maintheorem3.eq3}
		\|{\bf y}_{1, N}(t)-{\bf x}(t)\|_\infty \le \frac{\epsilon_0 \max \left\{ \|{\bf x}_0\|_2, R\epsilon_1 \right\}}{\eta_1(1-\epsilon_0)}
 \left(\frac{ \max\left\{\|{\bf x}_0\|_2, R\epsilon_1\right\}}{R\epsilon_0} 
  \right)^{N}
	\end{equation}
hold for all $t\ge t_0$ and $N\ge 1$, where ${\bf y}_{1, N}, N\ge 1$ is the first block in the finite section approximation \eqref{Athomogenous.def},
and
\begin{equation}\label{maintheorem3.eq3+}
\epsilon_0:= \frac{ \eta_0+\eta_1+\sqrt{  (\eta_1-\eta_0)^2- 4\eta_0  }} { 2(1+\eta_1)}\le \frac{\eta_1}{1+\eta_1}  \ {\rm and} \ \
\epsilon_1:=\frac{ \eta_0}{ (1+\eta_1) \epsilon_0}<\epsilon_0.
\end{equation}
\end{theorem}

In
\cite[Lemma 1 in the Supplementary Information]{liu2021efficient}, the authors consider
the Carleman linearization of nonlinear system  \eqref{dynamicsystem}
with ${\bf f}(t, {\bf x})$  being a  
polynomial
${\bf p}_{L}(t, {\bf x})$ as in \eqref{maincorollary.eq1} with $L=2$
and show that,
under similar assumptions
on the perturbation term ${\bf f}(t, {\bf 0})$  and  the initial  ${\bf x}_0$ to
the ones in \eqref{maintheorem3.eq1} and \eqref{maintheorem3.eq2},
there exists a positive constant $C$ such that
$\|{\bf y}_{1, N}(t)-{\bf x}(t)\|_2\le   C  t N  \|{\bf x}_0\|_2^{N}
$ hold
for all $t\ge t_0$ and $N\ge 1$.

For   nonlinear system \eqref{dynamicsystem} with the equilibrium assumption
${\bf f}(t, {\bf 0})= {\bf 0}$, we have
$\nu_0=\eta_0=\epsilon_1=0$ and
$\epsilon_0={R\mu_0}/({D_0+R\mu_0})$.
Hence,  the requirement \eqref{maintheorem3.eq1} is satisfied and
 the condition \eqref {maintheorem3.eq2} on  the initial  ${\bf x}_0$ boils down to
 \eqref{maintheorem2.eq1},
 and the  conclusion  \eqref{maintheorem3.eq3+} on the exponential convergence turns out to be  the same as the one in \eqref{maintheorem2.eq2}.
 Therefore,  Theorem \ref{maintheorem3} is a stronger version of Theorem \ref{maintheorem2}. Moreover, we should  emphasize that the exponential convergence in Theorem \ref{maintheorem3} does not imply the stability of  system
 \eqref{dynamicsystem} without the equilibrium assumption
${\bf f}(t, {\bf 0})= {\bf 0}$, i.e., the solution $\x(t)$ may not converge  as $t\to \infty$; however, it is shown in  Lemma \ref{maintheorem.lem2} that it is always bounded  by $\max\{\|{\bf x}_0\|, R\epsilon_1\}$.

\begin{remark} \label{maccoefficient.remark}
	An alternative hypothesis to  Assumption \ref{assump-1}
	is
	\begin{equation} \label{maccoefficient.remark.eq1}
		|f_{j, \pmb \alpha}(t)|\le \tilde D \tilde R^{-|\pmb \alpha|}
	\end{equation}
for all $1\le j\le d$,  ${\pmb \alpha}\in   {\mathbb Z}^d_+\backslash \{0\}$, and  $t\ge t_0$, 	where  $\tilde D$ and $\tilde R$  
	are some positive constants.
	Clearly, if 
	Assumption \ref{assump-1} holds, then the requirement \eqref{maccoefficient.remark.eq1} is satisfied with  $\tilde R=R$ and $\tilde D=D_0$
	using \eqref{fassump.eq1}. 
	Conversely, if  \eqref{maccoefficient.remark.eq1} is satisfied,
	then for any $ R<\tilde R$ there exists a positive constant $D_0$ such that
	the uniform exponential decay property \eqref{fassump.eq1} and,  hence,  Assumption \ref{assump-1} hold because of
	\begin{equation*}
		\sup_{t\ge t_0} \sum_{j=1}^d \sum_{\pmb \alpha\in { \mathbb Z}_k^d} |f_{j, \pmb\alpha}(t)|
		\le  d \tilde D {{k+d-1} \choose {d-1}}  \tilde R^{-k} \ {\rm and}
		\ \sup_{k\ge 1}   {{k+d-1} \choose {d-1}} \left(\frac{R}{\tilde R}\right)^{k}<\infty.
	\end{equation*}
	From the above observation and the conclusions in Theorems \ref{maintheorem}, \ref{maintheorem2} and \ref{maintheorem3},  under the uniform
	exponential decay assumption \eqref{maccoefficient.remark.eq1},  instead of using
	Assumption \ref{assump-1},
	we conclude that
	the first block  ${\bf y}_{1, N}$  of
	the solution of the finite-section approximation  \eqref{finitesection-1}  converges exponentially to the solution  ${\bf x}$ of the original nonlinear dynamical system
	\eqref{dynamicsystem}  over a quantfiable time interval  when the initial ${\bf x}(t_0)={\bf x}_0$ is in a vicinity of the equilibrium. Similarly, the exponential convergence will hold over the entire time horizon from $t_0$ to infinity  if \eqref{maccoefficient.remark.eq1} and Assumption \ref{assump-2} are both satisfied.
\end{remark}

\begin{remark}
In Assumption \ref{assump-2}, the requirement on the Jacobian matrix to be diagonal can be relaxed by requiring diagonalizability and  the similar results  can be established
with different constants in \eqref{maintheorem3.lem1.eq2}-\eqref{maintheorem3.lem2.eq1}. To circumvent this, one may first transform (rotate) the state variable of the original nonlinear system and rewrite  the corresponding dynamics in the new coordinates with diagonal Jacobian and then verify Assumption \ref{assump-2}.
\end{remark}

\section{Numerical simulations}
\label{numericalsimulation.section}
In this section, we study two nonlinear 
systems to validate and illustrate  conclusions
and  requirements
in Theorems \ref{maintheorem} and \ref{maintheorem2}. 

\subsection{Carleman linearization of  one-dimensional  dynamical systems}\label{unstable.simulation}
Let us consider  the following  one-dimensional nonlinear  dynamical systems
\begin{equation}
	\label{unstable.eq1}{\dot{x}}^\pm(t)= \pm \frac{ x^\pm (t)}{1+(x^\pm (t))^2}
\end{equation}
with the initial $x^\pm(0)=x_0>0$. One may verify that the solution of  the above dynamical system
is given by
\begin{equation}  \label{unstable.eq2}
	\ln |x^\pm(t)|+\frac{(x^\pm  (t))^2}{2} \mp t= \ln x_0+\frac{x_0^2}{2},\end{equation}
which implies that
$\lim_{t\to +\infty} x^+(t) t^{-1/2}=\sqrt{2}$ and $\lim_{t\to +\infty} x^-(t) e^t= x_0 \exp(x_0^2/2)$.
Therefore, the  dynamical system  \eqref{unstable.eq1} with positive sign is unstable, while the one with negative sign is stable.

\begin{figure}[t]
	\subfloat
	{\includegraphics[width = 2.4in]{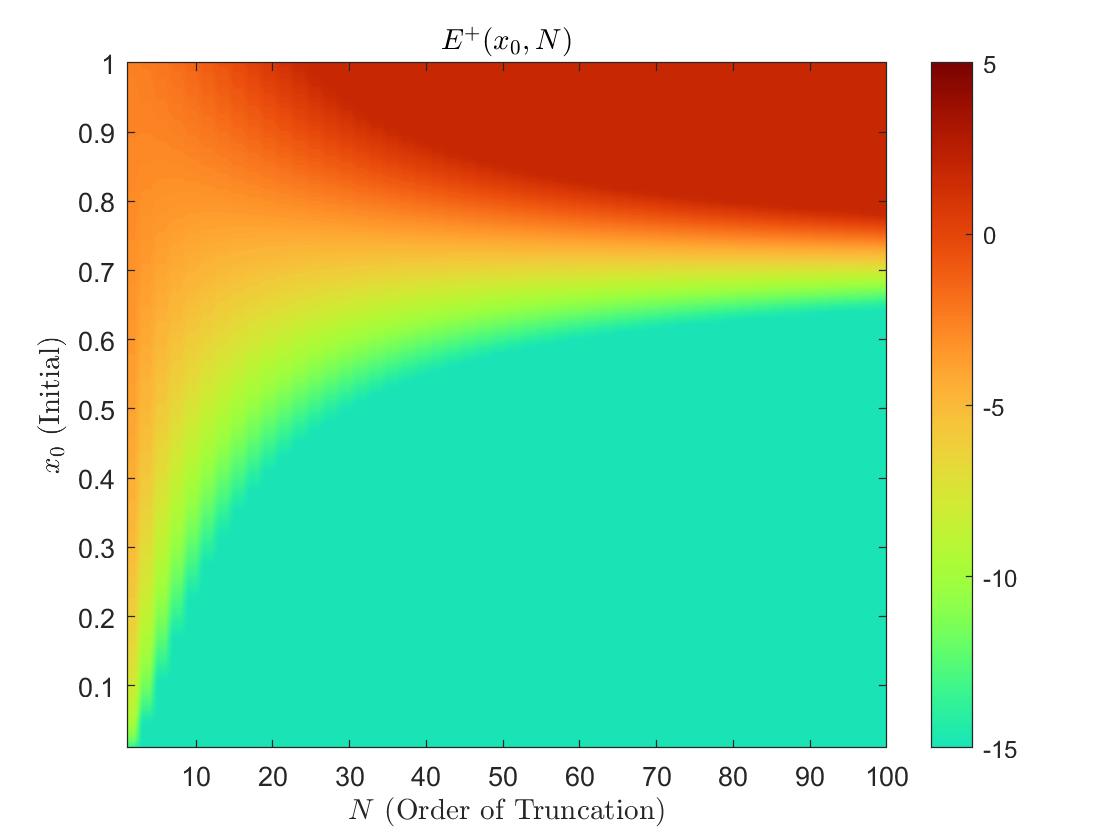}}
	\subfloat 
	{\includegraphics[width = 2.4in]{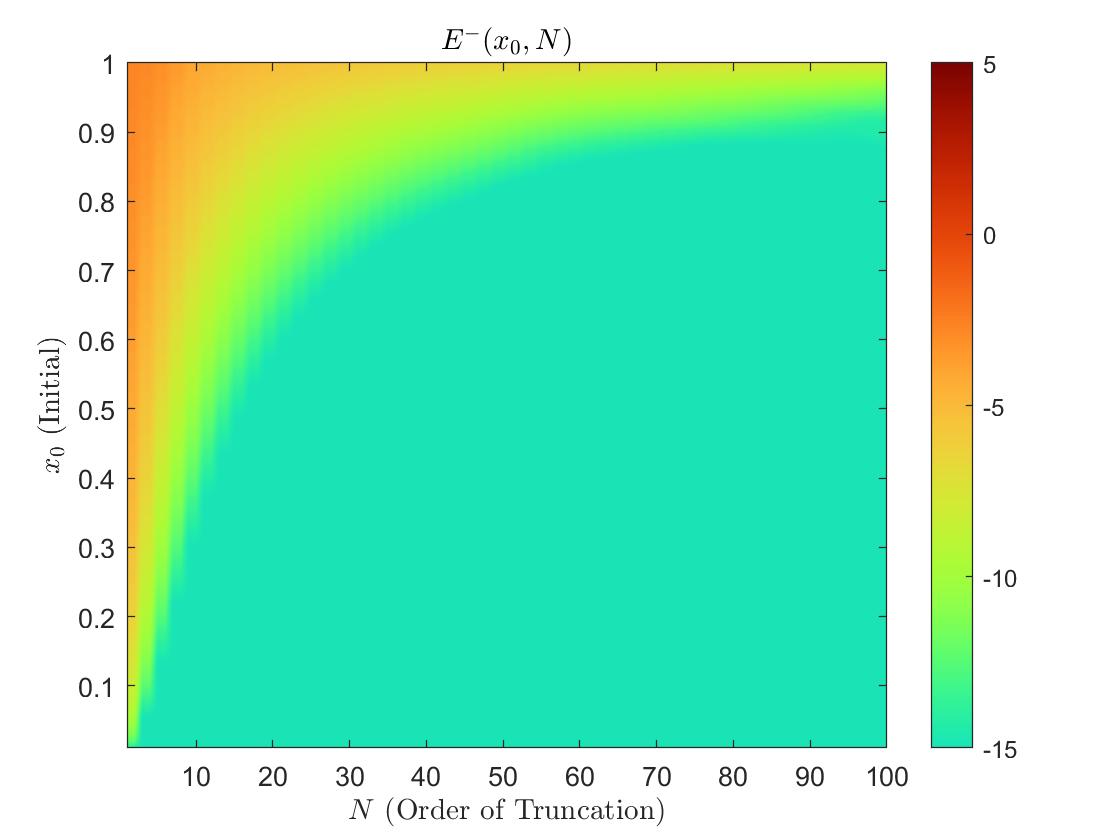}}
	\\
	\caption{Plotted  are maximal approximation error
		$E^\pm (x_0, N), 
		1\le N\le 100, 0\le x_0\le 1$
		in the logarithmic scale between
		the first block $y_{1, N}^\pm $ of the finite section scheme \eqref{unstable.eq4} 
		and the true solution $x^\pm$ of the  dynamical system
		\eqref{unstable.eq1} with positive sign (left) and   negative sign (right)  during the time period $[0, T^*_\pm]$.
	}
	\label{unstable.fig}
\end{figure}

By the Maclaurin expansion
${x}({1+x^2})^{-1}= \sum_{n=0}^\infty  (-1)^n x^{2n+1}$,
one can verify that  system  \eqref{unstable.eq1}, with respect to both signs, satisfies	 Assumption \ref{assump-1}  
with   $D_0=R=1$ and the dynamical system  \eqref{unstable.eq1} with negative sign also satisfies
Assumption \ref{assump-2} with $\mu_0=1$.
The  corresponding  finite-section scheme   of  the Carleman linearization of the  system \eqref{unstable.eq1}
is given by
\begin{equation}
	\label{unstable.eq4}
	\left [\begin{array}{c}
		\dot{ y}_{1, N}^\pm (t)\\
		\dot{ y}_{2, N}^\pm(t)\\
		\vdots\\
		\dot{ y}_{N, N}^\pm(t)
	\end{array}\right ]	 =\pm
	\left [\begin{array}{ccccc}
		1 &  0 & -1 &  \cdots  & ~ \cos \frac{(N-1)\pi}{2} \\
		0  & 2 & 0 &  \cdots &  2\cos \frac{(N-2)\pi}{2} \\
		0  & 0 & 3 &  \cdots &  3 \cos \frac{(N-3)\pi}{2} \\
		\vdots   & \vdots & \vdots & \ddots & \vdots \\
		0 &  0 & 0 & \cdots & N \cos \frac{(N-N)\pi}{2}
	\end{array}\right ]
	\left [\begin{array}{c}
		{y}_{1, N}^\pm(t)\\
		{ y}_{2, N}^\pm(t)\\
		\vdots\\
		{ y}_{N, N}^\pm(t)
	\end{array}\right ]  
\end{equation}
with $y_{k, N}^\pm(0)=x_0^k$ for  $1\le k\le N$.
For a given initial condition $x_0>0$, the truncation order $N\ge 1$, and time range $T^*>0$, let us denote the maximal (worst) approximation error between $y_{1, N}^\pm$ and  the true solution $x^\pm$
during the time interval $[0,  T^*]$ by

\[
e_\pm(x_0, N)= \sup_{0\le t\le T^*_\pm}
\left | y_{1, N}^\pm (t)-x^\pm (t) \right |,
\]
where \begin{equation*}  \label{maximalerror.def2}
	T^*_+=\max\left\{\frac{e-1}{2e-1}\left(\ln |x_0|^{-1}-1\right),  0.1\right\}\ \ {\rm and} \ \ T^*_-= 10.\end{equation*}
 Figure \ref{unstable.fig} shows
the maximal approximation error
$$  E^\pm (x_0, N)  =   \log \min \left\{ \max \left\{ e_\pm(x_0, N), 10^{-15}\right\}, 10^5 \right\}
$$
in the logarithmic scale.
We observe that
for  the  stable dynamical system  \eqref{unstable.eq1} with negative sign,
the state variable $y_{1, N}^-$ for truncation orders $1\le N\le 100$
approximates the true solution $x^-$ quite well during the 10 seconds time interval for all  initials
 $0<x_0<1$,
while for the unstable dynamical system  \eqref{unstable.eq1} with positive sign,
the state variable $y_{1, N}^+$ for truncation orders $1\le N\le 100$  approximates the true solution $x^+$ quickly during time interval $[0,T^*_+]$  for all initials $0<x_0<0.7$.
This validates the conclusions in Theorems \ref{maintheorem} and \ref{maintheorem2}
about exponential convergence
when the initial  is not far away from the origin.
As expected from the requirements
\eqref{maintheorem.eq2} and \eqref{maintheorem2.eq1} in Theorems \ref{maintheorem} and \ref{maintheorem2},
the state variable $y^+_{1, N}, 1\le N\le 100$,
does not provide good approximation to the true solution  $x^+$  even over time interval $[0,0.1]$ when the initial is  far away from the origin, e.g., when $x_0>0.8$.

\begin{figure}[t]
	\subfloat 
	{\includegraphics [width = 2.4in] 
		{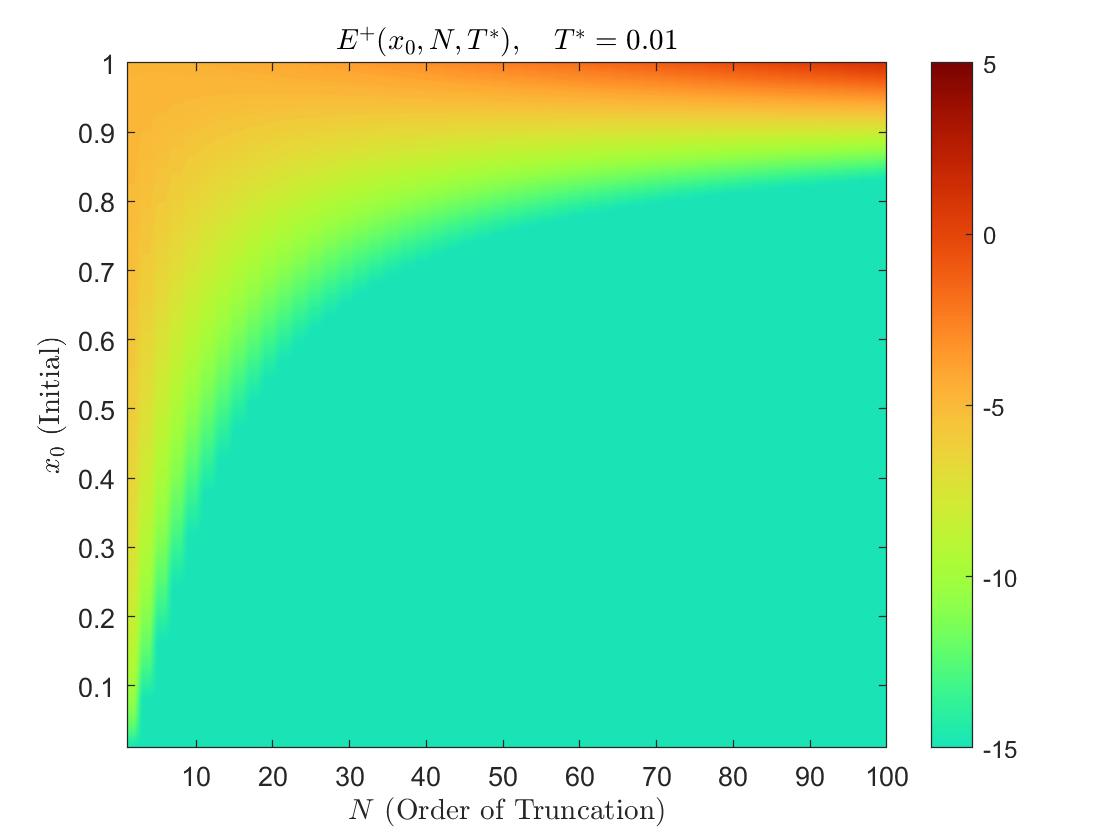}}
	\subfloat 
	{\includegraphics  [width = 2.4in]
		{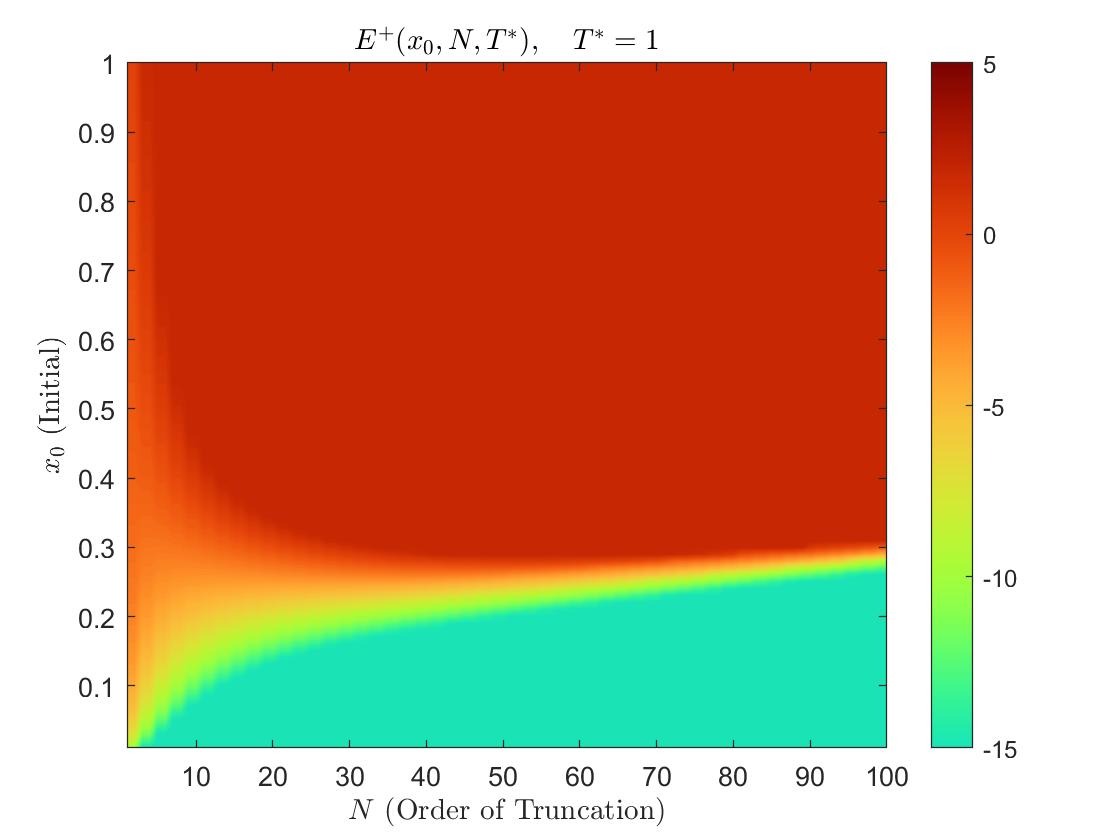}}\\
	\subfloat 
	{\includegraphics [width = 2.4in]{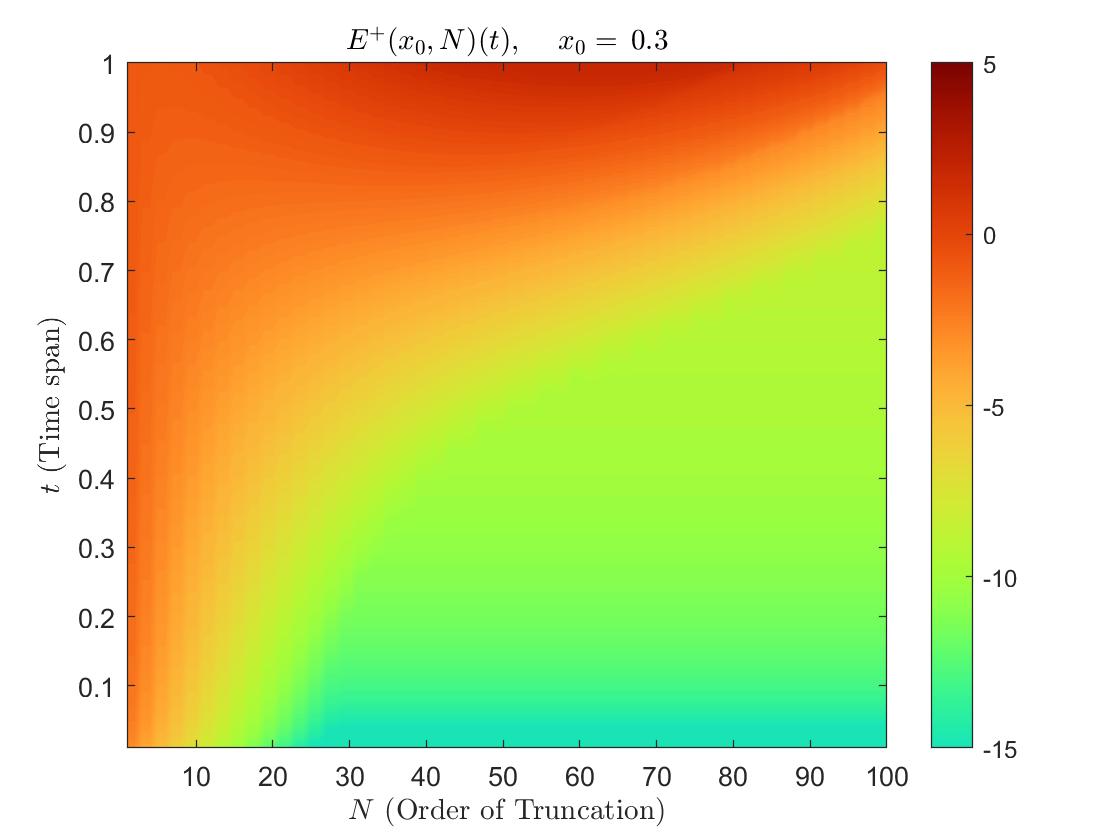}}
	\subfloat 
	{\includegraphics [width = 2.4in] {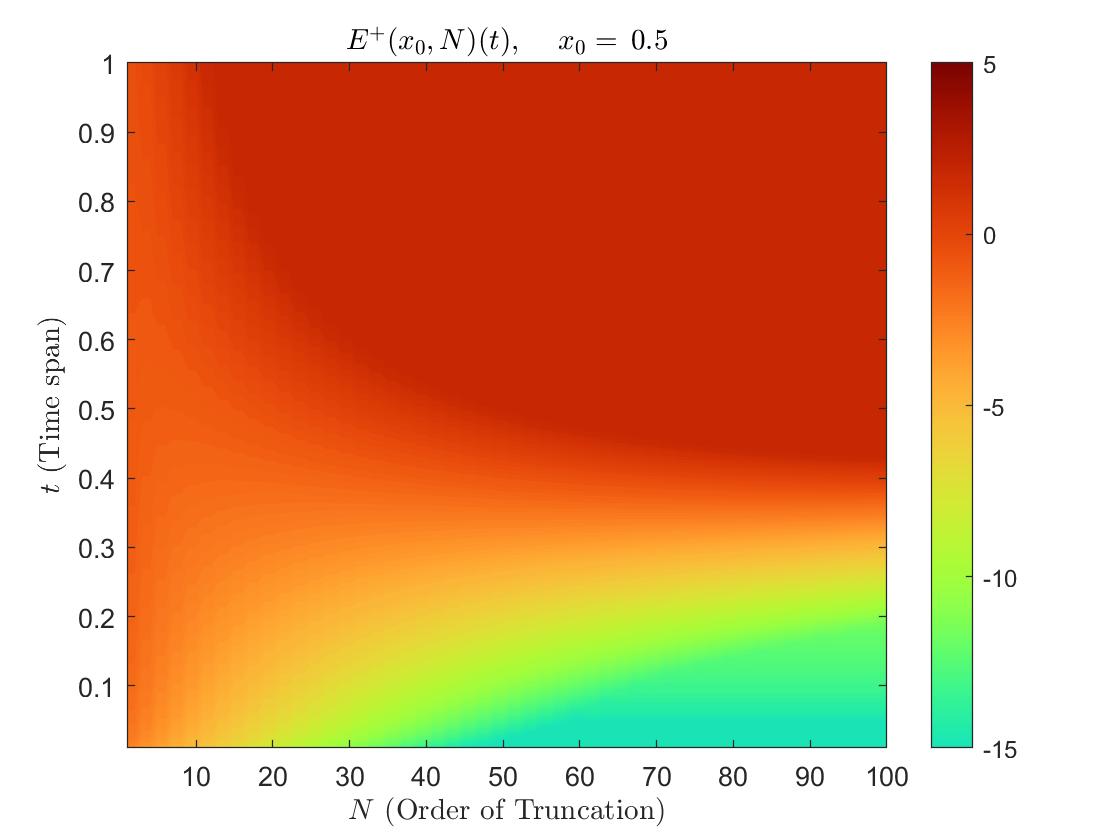}}
	\caption{Plotted on the top are maximal approximation error
		$E^+(x_0, N, T^*),  1\le N\le 100, 0\le x_0\le 1,$
		in the logarithmic scale
		with $T^*$ taking value $0.01$ (top left) and $1$ (top right) respectively.
		Plotted on the bottom  are  approximation error
		$E^+(x_0, N)(t), 1\le N\le 100, 0\le t\le 1$
		in the logarithmic scale for the initial $x_0=0.3$ and $0.5$ respectively.
	}
	\label{unstable.fig2}
\end{figure}

We  conduct further numerical simulations on the convergence rate for different time intervals $[0, T^*]$
and initials  $x_0$; see
Figure \ref{unstable.fig2}.
The top plots in Figure \ref{unstable.fig2}  are
maximal approximation error
$$E^+(x_0, N, T^*)=\log \min \left\{ \max  \left\{  e_+(x_0, N, T^*), 10^{-15}\right\}, 10^5 \right\}
$$
on time period $[0, T^*]$
in the logarithmic scale for  different time ranges $T^*$,
where the convergence region, i.e., $\log E^+(x_0,N, T^*)\le -15$, is marked by the green  color.
As expected from  \eqref{Tstar.def} in Theorem \ref{maintheorem},
increasing the convergence range $[0, T^*]$ results in smaller initial range $[0, x_0]$.
In particular, when  we increase $T^*$ from $0.01$ to $0.10$ and then to $1.00$, the maximal initial condition $x_0$ decreases
from $0.8423$ to $0.6722$ and then to $0.2348$,  while the theoretical bound in \eqref{Tstar.def} for the initial condition $x_0$ are $0.3582, 0.2841$, and $0.0401$, respectively.
The plots in the bottom of Figure \ref{unstable.fig2}
illustrate  the approximation error
$$E^+(x_0, N, t)=\log \min \left\{ \max \left\{ \left|y_{1, N}^+(t)-x^+(t)\right|, 10^{-15}\right\}, 10^5 \right\}
$$
for  $\ 0\le t\le 1$, in the logarithmic scale for various initials $x_0$.
 We notice that increasing the initial  $x_0$ leads to the shrinkage of the convergence region marked by the green color, or equivalently, the  time range for the convergence
gets smaller  when large  initials are selected. This asserts that the finite-section scheme \eqref{unstable.eq4} of the Carleman linearization should not be utilized to approximate the original
nonlinear dynamical system \eqref{unstable.eq1} with positive sign over long time intervals if the original nonlinear system is unstable.

\subsection{Carleman linearization of Van der Pol oscillator}\label{polynomial.simulation}	
\begin{figure}[t]
	\subfloat{\includegraphics[width =2.4in]{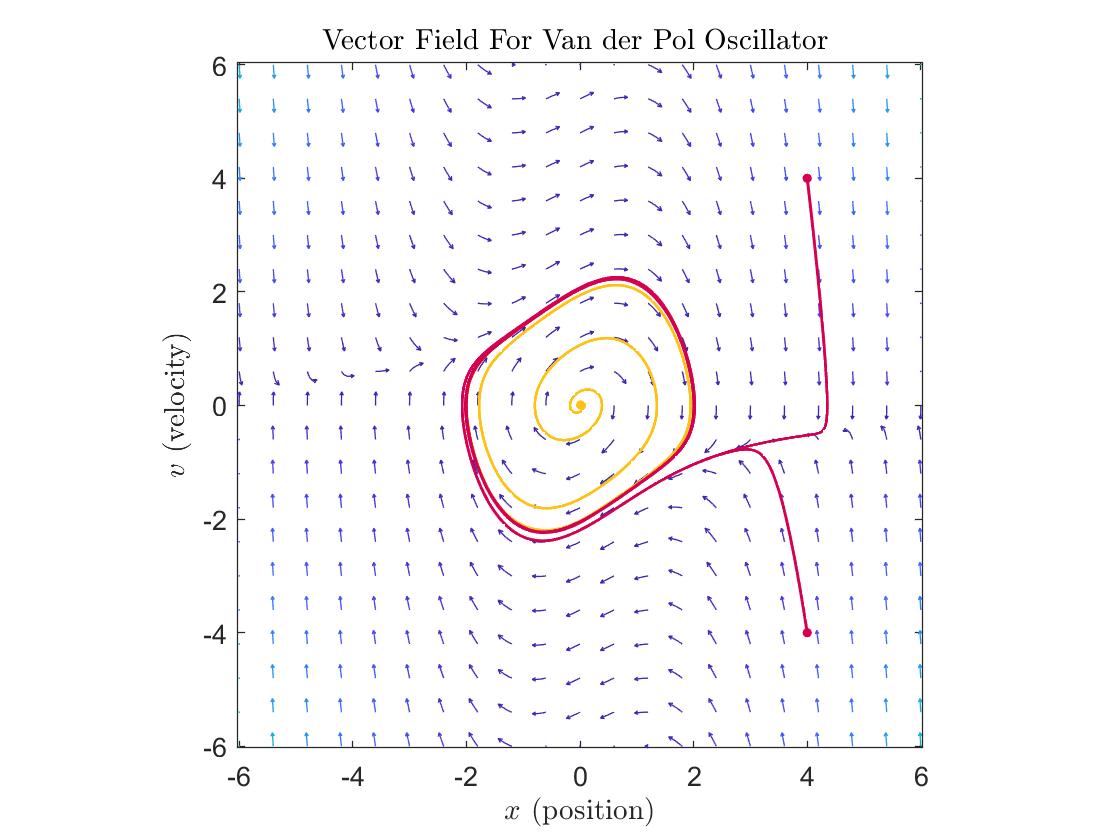}}
	\subfloat{\includegraphics[width = 2.4in]{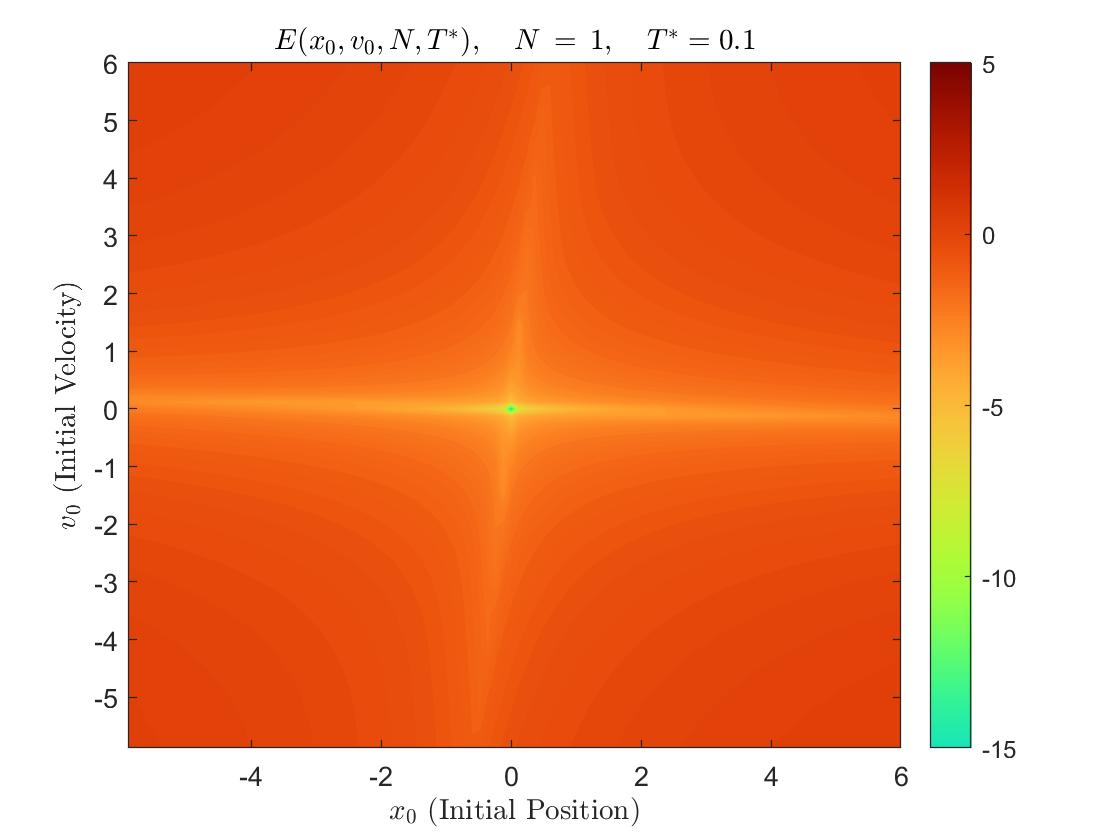}}\\
	\subfloat{\includegraphics[width = 2.4in] {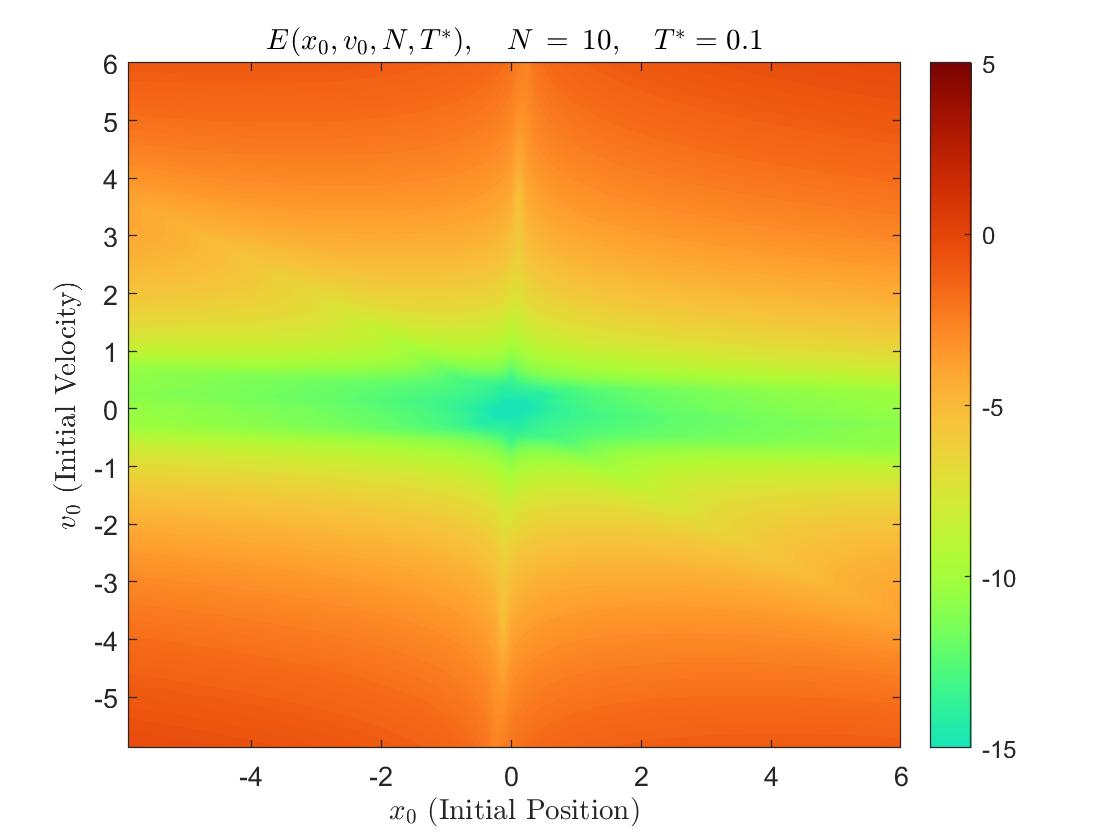}}
	\subfloat{\includegraphics[width = 2.4in] {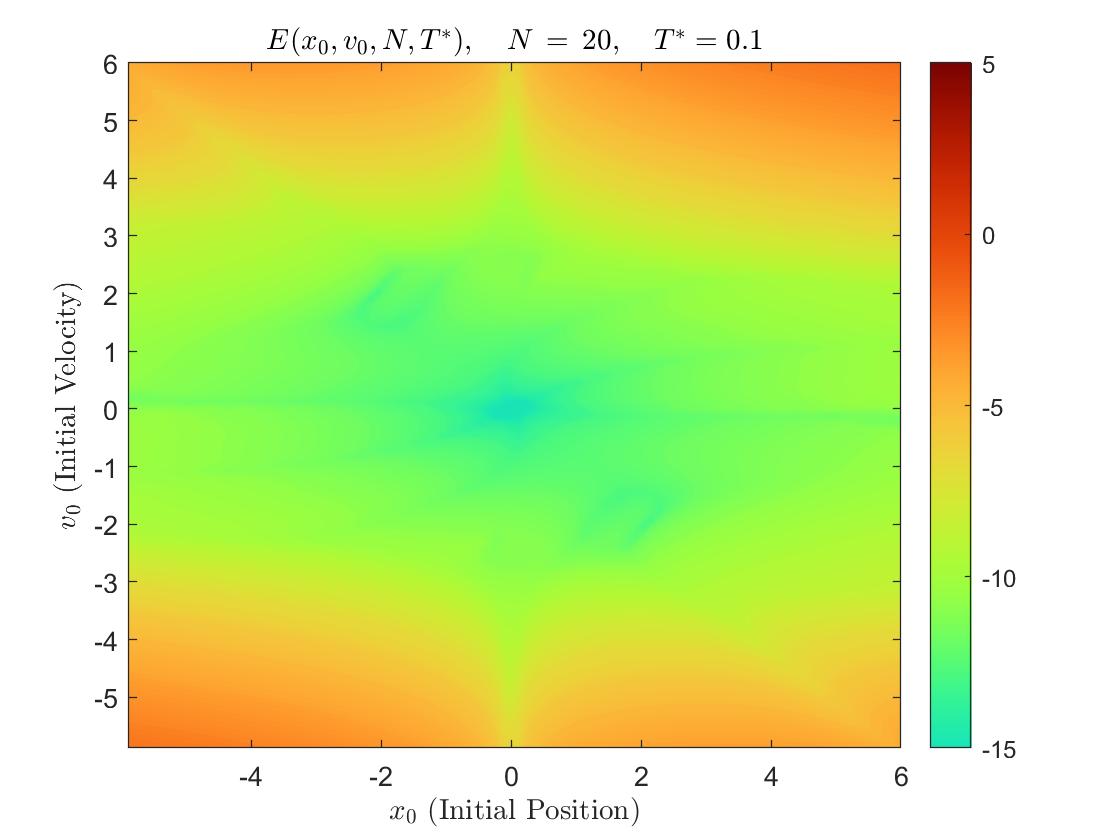}}	
	\caption{Plotted on the top left is  the
		vector field and limit cycle representation of the Van der Pol oscillator \eqref{vanderpol.def} with  $\mu=0.5$. On the top right, bottom left and  bottom right
		are the maximal approximation error $E(x_0,v_0,N,T^*)$
		in the logarithmic scale, where $T^*=0.1, \mu=0.5$  
		and
		the truncation order $N=1, 10, 20$  respectively.  }
	\label{unstable2.fig2}
\end{figure}

In this part, 
we consider Carleman linearization of
the
Van der Pol oscillator whose dynamics can be described by a dynamical system of polynomial type
\begin{equation}\label{vanderpol.def0}
	\ddot{x} -\mu(1-x^2)\dot{x} +x =0,
\end{equation}
where $x$ stands for position and $\mu$ is an indicator parameter for  non-linearity of the model.
Using the state vector
${\bf x}=[x_1, x_2]^T=[x, \dot{x}]^T$, this second-order  system can be represented in the canonical form \eqref{dynamicsystem} as
\begin{equation} \label{vanderpol.def}
	\dot{\bf x}=
	\left [ \begin{array}{c} x_2 \\
		-x_1+\mu(1-x_1^2) x_2\end{array}\right ]
	=
	\left [ \begin{array}{c} 0 \\
		-1\end{array}\right ] {\bf x}_{{\bf e}_1}
	+ \left [ \begin{array}{c} 1\\
		\mu\end{array}\right ] {\bf x}_{{\bf e}_2}
	+ \left [ \begin{array}{c} 0\\
		-\mu\end{array}\right] {\bf x}_{2{\bf e}_1+{\bf e_2}},
\end{equation}
where ${\bf e}_1=[1, 0]^T$ and ${\bf e}_2=[0, 1]^T$; see Figure \ref{unstable2.fig2} for the corresponding vector field. One may verify that the Van der Pol oscillator \eqref{vanderpol.def} satisfies Assumption \ref{assump-1} with 
$ d=2, R>0$, and $D_0=\max \left\{(2+\mu) R, \mu R^3 \right\}$.
By utilizing the new state variables  ${\bf z}_k=[x_1^k, x_1^{k-1}x_2, \ldots, x_2^k]^T$ for every  $k \ge 1$ and looking at the  sub-blocks of the state matrix \eqref{At.def}, it turns out that
${\bf A}_{k,l}$ for all $1\le k\le l$ will be zero matrices of size $(k+1)\times (l+1)$  except for
$${\bf A}_{k, k}= \left[\begin{array}{ccccccccc} 0 & k &   \\
	-1 & \mu & k-1 & \\
	& \ddots & \ddots & \ddots\\
	&  & \ddots & \ddots & \ddots\\
	& &  & \ddots & \ddots & \ddots\\
	&  & &  & -k+1& (k-1)\mu & 1 \\
	&  &  & & &  -k & k\mu
\end{array}
\right ] 
$$
and
$$
{\bf A}_{k,  k+2}= \left [\begin{array}{ccccccccc} 0&  &   & &0 & 0    \\
	&  -\mu &  & &  \vdots & \vdots \\
	&  & \ddots &  &\vdots& \vdots \\
	&  &  &  -k\mu & 0  & 0
\end{array}
\right]. 
$$
Therefore, the Carleman linearization of the  Van der Pol oscillator
is given by
\begin{equation}
	\left [\begin{array}{c}
		\dot {\bf z}_1\\
		\dot {\bf z}_2\\
		\vdots
	\end{array}\right]=
	\left[\begin{array}{cccccccccc}  {\bf A}_{1, 1} & {\bf 0}  & {\bf A}_{1, 3}  & \\
		& {\bf A}_{2, 2} & {\bf 0}  & {\bf A}_{2, 4}  & \\
		& & \ddots & \ddots & \ddots  & \\
	\end{array}\right]  \left [ \begin{array}{c}
		{\bf z}_1\\
		{\bf z}_2\\
		\vdots
	\end{array}\right ]
\end{equation}
and its corresponding  finite-section scheme of order $N$ is given by
\begin{equation}\label{VanderPol.finitesection}
	\left[\begin{array}{c}
		\dot {\bf y}_{1, N}\\
		\dot {\bf y}_{2, N}\\
		\vdots\\
		\vdots\\
		\dot {\bf y}_{N-1, N}\\
		\dot{\bf y}_{N, N}\\
	\end{array}\right]=
	\left[\begin{array}{cccccccccc}  {\bf A}_{1, 1} & {\bf 0}  & {\bf A}_{1, 3}  & \\
		& {\bf A}_{2, 2} & {\bf 0}  & {\bf A}_{2, 4}  & \\
		& & \ddots & \ddots & \ddots  & \\
		& & & \ddots & \ddots &  \\
		&  & & &   {\bf A}_{N-1, N-1} & {\bf 0} \\
		&  & & &    & {\bf A}_{N, N}  \\
	\end{array}\right]  \left [\begin{array}{c}
		{\bf y}_{1, N}\\
		{\bf y}_{2, N}\\
		\vdots\\
		\vdots\\
		{\bf y}_{N-1, N}\\
		{\bf y}_{N, N}
	\end{array}\right]
\end{equation}
with initial ${\bf y}_{k, N}(0)=[x_0^k, x_0^{k-1} v_0, \ldots, x_0 v_0^{k-1}, v_0^k]^T$,
where $x_0$ and $v_0$ are the initial position and velocity at time 0.

Let us define the maximal (worst) approximation error of the finite section  scheme \eqref{VanderPol.finitesection}
on the time interval $[0, T^*]$ in the logarithmic scale  by
$$E(x_0, v_0, N, T^*)=
\log \min \left\{ \max \left\{ \sup_{0\le t\le T^*} \left\|{\bf y}_{1, N}(t)- [x(t), x'(t)]^T \right\|_\infty, 10^{-15}\right\}, 10^5 \right\} $$
for $N\ge 1$,  where $ x (t)$ is the solution of
the Van der Pol oscillator \eqref{vanderpol.def0} with initial position $x_0$ and velocity $v_0$.
We verify the approximation properties of
the first block ${\bf y}_{1, N}(t)$  of the finite-section approximation  \eqref{VanderPol.finitesection}
with respect to the true solution $[x(t), \dot{x}(t)]^T$ over a fixed time period $[0,T^*]$ for different values of the truncation order $N$.
As expected from the exponential convergence
results in Theorem \ref{maintheorem},
the convergence region of the  initial $[x_0, v_0]^T$
is larger for  bigger truncation orders. Figure \ref{unstable2.fig2} shows
the maximal approximation error  $E(x_0, v_0, N, T^*)$ for the range of initial conditions  $-6\le x_0, v_0\le 6$.
As observed from the simulations,
the convergence region for the  traditional first-order linearization approach, i.e., when $N=1$,
is much  smaller (almost negligible) than the one for the finite-section approximation with large truncation orders, e.g.,  $N=10$ or $20$, which is depicted in Figure \ref{unstable2.fig2}.
These simulations assert that
Carleman linearization of nonlinear system \eqref{dynamicsystem}  can tightly approximate the original system over larger
regions and longer time intervals with respect to the conventional first-order linearization approach.

\section{Technical proofs of the main theorems} 
\label{proof.section}
In this section, we prove Theorems
\ref{maintheorem}
and \ref{maintheorem3}.

\subsection{Proof of Theorem \ref{maintheorem}} \label{maintheorem.pfsection}
For two countable index sets  $X$ and $Y$, let
$\ell^p(X)$  be the space of all $p$-summable sequences ${\bf u}=[u(i)]_{i\in X}$ with norm denoted  by $\|{\bf u}\|_{\ell^p(X)}$,
and  ${\mathcal S}(X, Y)$ and
${\mathcal B}_p(X, Y), 1\le p\le \infty$, be  Banach spaces of matrices  ${\bf C}= [C(i,j)]_{i\in X, j\in Y}$
with the norm defined by
\begin{equation}
	\label{maximalSchur.def}
	\|{\bf C}\|_{ {\mathcal S}(X, Y)}:=\max\left\{ \sup_{i\in X} \sum_{j\in Y} |C(i,j)|, \
	\sup_{j\in Y} \sum_{i\in X} |C(i,j)|\right\}
\end{equation}
and
\begin{equation}\label{maximaloperatornorm.def}
	\|{\bf C}\|_{{\mathcal B}_p(X, Y)}
	:= \sup_{{\bf u}\ne {\bf 0}}\frac{\|{\bf C}{\bf u}\|_{\ell^p(X)}}{\|{\bf u}\|_{\ell^p(Y)}},
\end{equation}
respectively. The norms in \eqref{maximalSchur.def} and \eqref{maximaloperatornorm.def} are known as the Schur norm and
the operator norm from $\ell^p(Y)$ to $\ell^p(X)$.
By direct computation, one may verify that
$\|{\bf C}\|_{{\mathcal B}_1(X,Y)}=
\sup_{j\in Y} \sum_{i\in X} |C(i,j)|$ and $\|{\bf C}\|_{{\mathcal B}_\infty(X, Y)}=
\sup_{i\in X} \sum_{j\in Y} |C(i,j)|$.
This together with the  interpolation theory (\cite{bergh2012interpolation}) yields
\begin{equation}\label{Schur.eq1}
	\|{\bf C}\|_{{\mathcal S}(X,Y)}= \sup_{1\le p\le \infty} \|{\bf C}\|_{{\mathcal B}_p(X,Y)}.\end{equation}
For  countable index sets $X, Y, Z$ and matrices ${\bf C}\in {\mathcal S}(X,Y)$ and ${\bf D}\in {\mathcal S}(Y, Z)$, one has
\begin{eqnarray}\label{Schur.eq2}
	\|{\bf C}{\bf D}\|_{{\mathcal S}(X,Z)} & \hskip-0.01in = & \hskip-0.01in
	\sup_{1\le p\le \infty} \|{\bf C}{\bf D}\|_{{\mathcal B}_p(X,Z)}\nonumber\\
	& \hskip-0.01in \le & \hskip-0.01in \sup_{1\le p\le \infty} \|{\bf C}\|_{{\mathcal B}_p(X, Y)}
	\|{\bf D}\|_{{\mathcal B}_p(Y,Z)}
	\le \|{\bf C}\|_{{\mathcal S}(X,Y)}  \|{\bf D}\|_{{\mathcal S}(Y,Z)}
	.\end{eqnarray}
Then,  for any countable set $X$,  ${\mathcal S}(X, X)$ is a Banach subalgebra of ${\mathcal B}_p(X, X), 1\le p\le \infty$ 
which is known as the Schur algebra, see  \cite{grochenig2006symmetry,  motee2017sparsity, shin2019polynomial,  sun2005wiener,
	sun2007wiener} for the inverse-closed property of Banach algebra of infinite matrices.

To prove Theorem \ref{maintheorem}, we need  the following Schur norm estimates for block matrices
${\bf A}_{k,l}(t)$  of  the time-varying state matrices ${\bf A}(t)$ in \eqref{At.def}.

\begin{lemma}\label{maintheorem.lem1} Let ${\bf f}(t, {\bf x})$ be as in Theorem \ref{maintheorem}, and
	${\bf A}_{k,l}(t), k, l\ge 1$ be  defined by \eqref{carleman.eq2}. Then
	\begin{equation}  \label{maintheorem.lem1.eq1}
		{\bf A}_{k,l}(t)=0\ \ {\rm if}\ \  1\le l<k,
	\end{equation}
	and
	\begin{equation} \label{maintheorem.lem1.eq2}
		\sup_{t\ge t_0} \|{\bf A}_{k,l}(t)\|_{{\mathcal S}( {\mathbb Z}^d_k, {\mathbb Z}^d_l)}\le  D_0 k R^{k-l-1} \ \ {\rm if}\  \   1\le k\le l.
	\end{equation}
\end{lemma}

\begin{proof}
	First we prove \eqref{maintheorem.lem1.eq1}.
	By \eqref{fj.negative} and \eqref{carleman.eq2}, the  proof of the conclusion \eqref{maintheorem.lem1.eq1} reduces to establishing
	${\pmb \beta}-{\pmb \alpha}+{\bf e}_j\not\in {\mathbb Z}_+^d\backslash \{0\}, \ 1\le j\le d$,
	for all ${\pmb \alpha}\in {\mathbb Z}^d_k$ and ${\pmb \beta}\in {\mathbb Z}_l^d$ with $l<k$.
	Suppose, on the contrary, that
	${\pmb \beta}-{\pmb \alpha}+{\bf e}_j\in {\mathbb Z}_+^d\backslash \{0\}$ for some $1\le j\le d$. Then
	$l+1=|{\pmb \beta}+{\bf e}_j|= |{\pmb \alpha}+ ({\pmb \beta}-{\pmb \alpha}+{\bf e}_j)|>|{\pmb \alpha}|=k$,
	which  contradicts to the assumption that $1\le l<k$.
	
	Next we prove \eqref{maintheorem.lem1.eq2}. For ${\pmb \alpha}\in {\mathbb Z}^d_+$, we write ${\pmb \alpha}=[\alpha_1, \ldots, \alpha_d]^T$.
	For  integers $k$ and $l$ with $k\le l$ and  $t\ge t_0$,   
	\begin{eqnarray} \label{maintheorem.lem1.pfeq2}
		& & \hskip-0.01in \max\left\{ \sup_{{\pmb \alpha}\in {\mathbb Z}_k^d} \sum_{{\pmb \beta}\in {\mathbb Z}_l^d} \left|\sum_{j=1}^d \alpha_j f_{j, \pmb\beta-\pmb\alpha+{\bf e}_j}(t)\right|,\
		\sup_{{\pmb \beta}\in {\mathbb Z}_l^d} \sum_{{\pmb \alpha}\in {\mathbb Z}_k^d}\left |\sum_{j=1}^d \alpha_j f_{j, \pmb\beta-\pmb\alpha+{\bf e}_j}(t)\right|\right\}\nonumber\\
		&  \hskip-0.01in \le &  \hskip-0.01in k   \max\left\{ \sup_{{\pmb \alpha}\in {\mathbb Z}_k^d} \sum_{{\pmb \beta}\in {\mathbb Z}_l^d} \sum_{j=1}^d \left|f_{j, \pmb\beta-\pmb\alpha+{\bf e}_j}(t)\right|,\
		\sup_{{\pmb \beta}\in {\mathbb Z}_l^d} \sum_{{\pmb \alpha}\in {\mathbb Z}_k^d} \sum_{j=1}^d |f_{j, \pmb\beta-\pmb\alpha+{\bf e}_j}(t)|\right\}\nonumber\\
		& \hskip-0.01in \le &  \hskip-0.01in  k
		\sum_{j=1}^d   \sum_{{\pmb \gamma}\in {\mathbb Z}^d_{l-k+1}} \left|f_{j, \pmb \gamma}(t)\right|
		\le  D_0 k R^{k-l-1},
	\end{eqnarray}
	where the first inequality holds as $0\le \alpha_j\le  |{\pmb \alpha}|=k$ for all $1\le j\le d$,
	the second inequality follows from  \eqref{fj.negative} and the observation that
	$|{\pmb \beta}-{\pmb \alpha}+{\bf e}_j|=|{\pmb \beta}|-|{\pmb \alpha}|+1=k-l+1$ when ${\pmb \beta}-{\pmb \alpha}+{\bf e}_j\in {\mathbb Z}_+^d$,
	and the last inequality is true by the assumption
	\eqref{fassump.eq1}.  Taking supernum over $t\ge t_0$ in \eqref{maintheorem.lem1.pfeq2} proves  \eqref{maintheorem.lem1.eq2}.
\end{proof}

To prove Theorem \ref{maintheorem}, we need  the following local bound estimate for the solution ${\bf x}(t)$ of the nonlinear dynamical system
\eqref{dynamicsystem}.

\begin{lemma}
	\label{maintheorem.lem2}
	Let  ${\bf f}(t, {\bf x})$ be as in Theorem \ref{maintheorem} and ${\bf x}(t), t\ge t_0$, be
	the solution of the nonlinear dynamical system
	\eqref{dynamicsystem}.
	If the nonzero initial ${\bf x}_0$ satisfies \eqref{maintheorem.eq2},
	then \eqref{maintheorem.lem2.eq1} holds.
\end{lemma}

\begin{proof}  Let
	\begin{equation}  \label{maintheorem.lem2.pfeq1}
		t_1^*=\sup\big\{ t,\  \|{\bf x}(s)\|_\infty\le M_0 \ \ {\rm  for \ all} \  \ t_0\le s\le t\big\}.
	\end{equation}
	By the continuity of the solution ${\bf x}(t)$ of the nonlinear  system  \eqref{dynamicsystem} and
	the assumption \eqref{maintheorem.eq2} on the initial ${\bf x}_0$ that
	$\|{\bf x}(t_0)\|_\infty=\|{\bf x}_0\|_\infty< M_0$, we have
	that
	$t_1^*>t_0$.
	
	The conclusion \eqref{maintheorem.lem2.eq1} is obvious if  $t_1^*=+\infty$. Then it remains to consider the case that  $t_0<t_1<\infty$.
	In this case, we obtain from the continuity of ${\bf x}(t), t\ge t_0$, that
	\begin{equation}  \label{maintheorem.lem2.pfeq2}
		\|{\bf x}(t_1^*)\|_\infty= M_0.
	\end{equation}
	Integrating the nonlinear dynamical system  \eqref{dynamicsystem} yields
	\begin{equation*}
		{\bf x}(t)={\bf x}_0+\int_{t_0}^t {\bf f}(s, {\bf x}(s)) ds={\bf x}_0+\int_{t_0}^t \sum_{\pmb \alpha\in  {\bf Z}_+^d\backslash \{0\}} f_{\pmb \alpha} (s) ({\bf x}(s))^{\pmb \alpha} ds,
	\end{equation*}
	where the last equality follows from the Maclaurin expansion \eqref{marcluurin.def} of the function ${\bf f}(t, {\bf x})$.  Hence for $t_0\le t\le t_1^*$, we have
	\begin{eqnarray}\label{maintheorem.lem2.pfeq3}
		\|{\bf x}(t)\|_\infty & \hskip-0.01in \le & \hskip-0.01in  \|{\bf x}_0\|_\infty+\int_{t_0}^t \sum_{k=1}^\infty \left(\sum_{j=1}^d \sum_{{\pmb \alpha}\in {\mathbb Z}_k^d} |f_{j, \pmb \alpha}(s)|\right) \|{\bf x}(s)\|_\infty^k ds\nonumber\\
		&  \hskip-0.01in  \le  &  \hskip-0.01in  \|{\bf x}_0\|_\infty+ D_0 \int_{t_0}^t \sum_{k=1}^\infty R^{-k} \|{\bf x}(s)\|_\infty^k ds 
		\le \|{\bf x}_0\|_\infty+\frac{ D_0 }{ R-M_0}  \int_{t_0}^t  \|{\bf x}(s)\| ds,
	\end{eqnarray}
	where the second estimate follows from \eqref{fassump.eq1} and the last inequality holds by   \eqref{maintheorem.lem2.pfeq1}. Let us set
	$$u(t)=  \int_{t_0}^t  \|{\bf x}(s)\|_\infty ds-\frac{(R-M_0)\|{\bf x}_0\|_\infty}{D_0}  \left( e^{D_0(t-t_0)/(R-M_0)}-1\right), \ \ t_0\le t\le t_1^*. $$
	Then,  $u(t_0)=0$
	and
	\begin{eqnarray*}  & & \hskip-0.01in  \frac{ d\big(e^{-D_0(t-t_0)/(R-M_0)} u(t)\big)}{dt}  =
		e^{-D_0(t-t_0)/(R-M_0)} \left( -\frac{D_0}{R-M_0} u(t)+ \frac{ d u(t)}{dt}\right)\nonumber\\
		& \hskip-0.01in  \le & \hskip-0.01in
		e^{-D_0(t-t_0)/(R-M_0)} \left( \|{\bf x}(t)\|_\infty - \|{\bf x}_0\|_\infty-\frac{D_0}{R-M_0} \int_{t_0}^t \|{\bf x}(s)\|_\infty ds\right)
		\le 0,
	\end{eqnarray*}
	where the last inequality follows from \eqref{maintheorem.lem2.pfeq3}.
	This implies that
	$	u(t)\le 0$ for all $ t_0\le t\le t_1^*$,
	and hence
	\begin{equation}\label{maintheorem.lem2.pfeq4}
		\int_{t_0}^t  \|{\bf x}(s)\|_\infty ds\le \frac{(R-M_0)\|{\bf x}_0\|_\infty}{D_0}  \left( e^{D_0(t-t_0)/(R-M_0)}-1\right) \ {\rm for \ all} \ t_0\le t\le t_1^*.
	\end{equation}
	Substituting the above estimate into \eqref{maintheorem.lem2.pfeq3}, we obtain
	\begin{equation}\label{maintheorem.lem2.pfeq5}
		\|{\bf x}(t)\|_\infty\le \|{\bf x}_0\|_\infty   e^{D_0(t-t_0)/(R-M_0)}, \ t_0\le t\le t_1^*.
	\end{equation}
	Combining \eqref{maintheorem.eq2}, \eqref{maintheorem.M0def}, \eqref{maintheorem.lem2.pfeq2} and \eqref{maintheorem.lem2.pfeq5} proves that
	$$ t_1^*\ge t_0+ \frac{R-M_0}{D_0} \ln \frac{M_0}{\|{\bf x}_0\|_\infty}\ge t_0+ T^*. $$
	This together with  \eqref{maintheorem.lem2.pfeq1} proves
	\eqref{maintheorem.lem2.eq1}.
\end{proof}

\smallskip
To prove Theorem \ref{maintheorem}, we need a kernel estimate for the solution of an ordinary differential equation with the bounded state matrix in the Schur norm.

\begin{lemma}\label{maintheorem.lem3}
	Let  $1\le p\le \infty$, $X$ be a countable index set,  the vector-valued function ${\bf w}(t), t\ge t_0$, be locally bounded in $\ell^p(X)$,
	\begin{equation} \label{maintheorem.lem3.eq0a}
		\sup_{t_0\le t\le t_1} \|{\bf w}(t)\|_{\ell^p(X)}<\infty \ {\rm for \ all} \ t_1<\infty,
	\end{equation}
	and let the matrix-valued function
	${\bf B}(t), t\ge t_0$, be   bounded in the Schur algebra ${\mathcal S}(X, X)$,
	\begin{equation} \label{maintheorem.lem3.eq0}
		\|{\bf B}\|_{\infty, {\mathcal S}(X, X)}:=\sup_{t\ge t_0} \|{\bf B}(t)\|_{{\mathcal S}(X, X)}<\infty.
	\end{equation}
	Then the locally bounded  solution of the  ordinary different equation
	\begin{equation}\label{maintheorem.lem3.eq1}
		\dot{\bf z}(t)= {\bf B}(t) {\bf z}(t)+ {\bf w}(t), \  t\ge t_0  \ \ {\rm with \ zero\  initial}\ \  {\bf z}(t_0)={\bf 0}
	\end{equation}
	in $\ell^p(X)$ has the following expression,
	\begin{equation}\label{maintheorem.lem3.eq2}
		{\bf z}(t)=\int_{t_0}^t {\bf K}(t, s) {\bf w}(s) ds
	\end{equation}
	where
	the integral kernel ${\bf K}$ satisfies
	\begin{equation}  \label{maintheorem.lem3.eq3}
		\|{\bf K}(t, s)\|_{{\mathcal S}(X, X)}\le  \exp\big(\|{\bf B}\|_{\infty, {\mathcal S}(X, X)}(t-s)\big)
		\ \ {\rm for \ all}\ \  t\ge s\ge t_0.
	\end{equation}
\end{lemma}

\begin{proof} Define integral kernels ${\bf T}_n$ and ${\bf S}_n, n\ge 0$, inductively by
	$$
	{\bf T}_{n}(t, s):= \int_s^t {\bf T}_{n-1}(t, u) {\bf B}(s)  du\ \ {\rm  and} \  \
	{\bf S}_n(t, s):=  \int_s^t  {\bf T}_{n-1}(t, u) du, \  n\ge 1,
	$$
	with the initial integral kernels ${\bf T}_0$ and ${\bf S}_0$ given by
	$   {\bf T}_0(t, s):=   {\bf  B}(s)$ and ${\bf S}_0(t, s):= {\bf I}, \ \ t\ge s\ge t_0.$
	By induction, we may verify that
	\begin{equation} \label{maintheorem.lem3.pfeq1}
		\|{\bf T}_n(t, s)\|_{{\mathcal S}(X, X)}\le \|{\bf B}\|_{\infty, {\mathcal S}(X, X)}^{n+1} \frac{ (t-s)^n}{n!},
	\end{equation}
	and
	\begin{equation}  \label{maintheorem.lem3.pfeq2}
		\|{\bf S}_n(t, s)\|_{{\mathcal S}(X, X)}\le \|{\bf B}\|_{\infty, {\mathcal S}(X, X)}^{n} \frac{ (t-s)^n}{n!}, \ t\ge s \ge t_0,
	\end{equation}
	for all $n\ge 0$.
	Define
	\begin{equation}   \label{maintheorem.lem3.pfeq2+}
		{\bf K}(t, s):=\sum_{k=0}^\infty {\bf S}_k(t, s),\ \  t\ge s \ge t_0.\end{equation}
	Then we get from \eqref{maintheorem.lem3.pfeq2} that
	\begin{equation} \label{maintheorem.lem3.pfeq2++}
		\|{\bf K}(t, s)\|_{{\mathcal S}(X, X)}\le  \sum_{n=0}^\infty \|{\bf S}_k(t, s)\|_{{\mathcal S}(X, X)}
		\sum_{k=0}^\infty \|{\bf B}\|_{\infty, {\mathcal S}(X, X)}^{k} \frac{ (t-s)^k}{k!} 
	\end{equation}
	for all $t\ge s\ge t_0$. This  proves \eqref{maintheorem.lem3.eq3}.

	Now it remains to prove \eqref{maintheorem.lem3.eq2} with the integral kernel ${\bf K}$   given in \eqref{maintheorem.lem3.pfeq2+}.
	Integrating both sides of the ordinary differential equation \eqref{maintheorem.lem3.eq1}
	gives
	\begin{equation}\label{maintheorem.lem3.pfeq3-}
		{\bf z}(t)=\int_{t_0}^t {\bf B}(s) {\bf z}(s) ds + \int_{t_0}^t  {\bf w}(s)ds, \ t\ge t_0.
	\end{equation}
	Applying \eqref{maintheorem.lem3.pfeq3-} for $n$ times, we obtain
	\begin{eqnarray}\label{maintheorem.lem3.pfeq3}
		{\bf z}(t) &\hskip-0.01in  = & \hskip-0.01in \int_{t_0}^t {\bf B}(u) \left(\int_{t_0}^u   {\bf B}(s) {\bf z}(s) ds+ \int_{t_0}^u {\bf w}(s) ds\right) du + \int_{t_0}^t  {\bf w}(s)ds\nonumber\\
		&\hskip-0.01in  = & \hskip-0.01in \int_{t_0}^t {\bf T}_1(t, s) {\bf z}(s) ds+
		\int_{t_0}^t \left(\sum_{k=0}^1 {\bf S}_k(t, s\right) {\bf w}(s) ds= \cdots\nonumber\\
		&\hskip-0.01in  = & \hskip-0.01in \int_{t_0}^t {\bf T}_{n}(t, s) {\bf z}(s) ds+  \int_{t_0}^t \left(\sum_{k=0}^n {\bf S}_k(t, s)\right) {\bf w}(s) ds, \ t\ge t_0.
	\end{eqnarray}
	Using the local boundedness of ${\bf z}(t), t\ge t_0$, in $\ell^p(X)$
	and  applying \eqref{Schur.eq1} and
	\eqref{maintheorem.lem3.pfeq1}, we have
	\begin{eqnarray} \label{maintheorem.lem3.pfeq4}
		& \hskip-0.01in  & \hskip-0.01in  \Big\|\int_{t_0}^t {\bf T}_{n}(t, s) {\bf z}(s) ds \Big\|_{\ell^p(X)}
		\le   \int_{t_0}^t \|{\bf T}_{n}(t, s)\|_{{\mathcal S}(X, X)} \|{\bf z}(s)\|_{\ell^p(X)} ds\nonumber\\
		& \hskip-0.01in \le & \hskip-0.01in
		\frac{ \|{\bf B}\|_{\infty, {\mathcal S}(X, X)}^{n+1}  (t-t_0)^n}{n!} \int_{t_0}^t \|{\bf z}(s)\|_{\ell^p(X)} ds
		\to 0\ \ {\rm as} \ \ n\to \infty.
	\end{eqnarray}
	Similarly,  from \eqref{Schur.eq1}, \eqref{maintheorem.lem3.eq0a} and
	\eqref{maintheorem.lem3.pfeq2} we obtain
	\begin{eqnarray} \label{maintheorem.lem3.pfeq5}
		& \hskip-0.01in  & \hskip-0.01in  \Big\|\int_{t_0}^t  {\bf K}(t, s) {\bf w}(s) ds - \int_{t_0}^t \left(\sum_{k=0}^n {\bf S}_k(t, s)\right) {\bf w}(s) ds \Big\|_{\ell^p(X)}\nonumber\\
		& \hskip-0.01in \le & \hskip-0.01in
		\int_{t_0}^t  \sum_{k=n+1}^\infty \|{\bf S}_{k}(t, s)\|_{{\mathcal S}(X, X)} \|{\bf w}(s)\|_{\ell^p(X)} ds\nonumber\\
		& \hskip-0.01in \le & \hskip-0.01in
		\left( \sum_{k=n+1}^\infty
		\frac{ \|{\bf B}\|_{\infty, {\mathcal S}(X, X)}^{k}  (t-t_0)^k}{k!}\right) \int_{t_0}^t \|{\bf w}(s)\|_{\ell^p(X)} ds
		\to 0\ \ {\rm as} \ \ n\to \infty.
	\end{eqnarray}
	Combining \eqref{maintheorem.lem3.pfeq3}, \eqref{maintheorem.lem3.pfeq4} and \eqref{maintheorem.lem3.pfeq5}
	proves
	\eqref{maintheorem.lem3.eq2}. This together with \eqref{maintheorem.lem3.pfeq2++} completes the proof.
\end{proof}

The final technical lemma used in our proof of Theorem \ref{maintheorem}
is as follows.

\begin{lemma} \label{maintheorem.lem4}
	Let  $D_1>0$ and
	$v_k, 1\le k\le N$, be nonnegative functions satisfying
	\begin{equation}\label{maintheorem.lem4.pf1}
		0\le v_k(t)\le   D_1 k \int_{t_0}^t e^{D_1 k (t-s)} \left(\sum_{l=k+1}^N   v_l(s) +1 \right) ds,\ \  t\ge t_0,
	\end{equation}
 then
	\begin{equation} \label{maintheorem.lem4.pf2}
		v_{N-n}(t)+\cdots+v_N(t)+1\le
		\frac{N^n}{n!}  e^{D_1N(t-t_0)}, \ t\ge t_0
	\end{equation}
	for all $ 0\le n\le N-2$.
\end{lemma}

\begin{proof} We prove \eqref{maintheorem.lem4.pf2}
	by induction on $ n=0, 1, \ldots, N-2$.
	Applying \eqref{maintheorem.lem4.pf1} with $k=N$, we have
	\begin{equation*} \label{maintheorem.pfeq6+} 
		v_N(t)+1\le D_1 N \int_{t_0}^t e^{D_1N (t-s)}  ds+1 =  e^{D_1N (t-t_0)},\  t\ge t_0.
	\end{equation*}
	This proves  \eqref{maintheorem.lem4.pf2} for $n=0$.
	
	Inductively we assume that the conclusion  \eqref{maintheorem.lem4.pf2} for  $0\le n\le N-2$.
	Applying \eqref{maintheorem.lem4.pf1} with $k$ replaced by
	$N-n-1$, we obtain from the inductive hypothesis  that
	\begin{eqnarray*}
		& & v_{N-n-1}(t) +v_{N-n}(t)+\cdots+v_N(t)+1\nonumber\\
		&\hskip-0.01in  \le &  \hskip-0.01in D_1 (N-n-1)  \int_{t_0}^t e^{D_1(N-n-1)(t-s)}  \frac{N^n}{n!} e^{D_1N (s-t_0)} ds+
		\frac{N^n}{n!} e^{D_1N (t-t_0)}
		\\
		&\hskip-0.01in = & \hskip-0.01in  \frac{N^n(N-n-1)}{(n+1)!}  
		\Big( e^{D_1 N(t-t_0)}-  e^{D_1 (N-n-1)(t-t_0)}\Big)+
		\frac{N^n}{n!} e^{D_1 N(t-t_0)}
		\\
		& \hskip-0.01in \le  & \hskip-0.01in
		 \frac{N^{n+1}}{(n+1)!}  e^{D_1 N (t-t_0)},
	\end{eqnarray*}
	 for all $t\ge t_0$.
	Therefore the inductive argument can proceed. This completes the inductive proof of the conclusion \eqref{maintheorem.lem4.pf2}.
\end{proof}

We finish this subsection with the proof of Theorem \ref{maintheorem}.

\begin{proof}
	Set $\pmb \eta_{k,N}(t)={\bf y}_{k,N}(t)-{\bf z}_k(t), t\ge t_0$, for $1\le k\le N$
	and
	$\|{\bf A}_{k, l}\|_{\infty,{\mathcal S}({\mathbb Z}_k^d, {\mathbb Z}_l^d)}=
	\sup_{t\ge t_0} \|{\bf A}_{k, l}(t)\|_{{\mathcal S}({\mathbb Z}_k^d, {\mathbb Z}_l^d)}$, $ k, l\ge 1$.
	Then we obtain from
	\eqref{At.def}  that
	\begin{equation}\label{maintheorem.pfeq1}
		\dot{ \pmb \eta}_{k,N}(t)= {\bf A}_{k, k}(t) \pmb \eta_{k, N}(t)+ \sum_{l=k+1}^N {\bf A}_{k, l}(t) \pmb \eta_{l, N}(t) -\sum_{l= N+1}^\infty {\bf A}_{k,l}(t) {\bf z}_l(t) \end{equation}
	with  zero initial ${\bf \eta}_{k,N}(t_0)={\bf 0}, 1\le k\le N$.
	Applying  Lemma \ref{maintheorem.lem3} with ${\bf B}(t)$ and ${\bf w}(t)$ replaced by
	${\bf A}_{k, k}(t)$ and $\sum_{l=k+1}^N {\bf A}_{k, l}(t) {\pmb \eta}_{l, N}(t) -\sum_{l= N+1}^\infty {\bf A}_{k,l}(t) {\bf z}_l(t)$ respectively,
	we can construct an integral kernel
	${\bf K}_k(t, s), t, s \ge t_0$ such that
	\begin{equation} \label{maintheorem.pfeq2}
		{\pmb \eta}_{k, N}(t)=\int_{t_0}^t  {\bf K}_k(t, s)  \left(\sum_{l=k+1}^N {\bf A}_{k, l}(s) {\pmb \eta}_{l, N}(s) -\sum_{l=N+1}^\infty {\bf A}_{k,l}(s) {\bf z}_l(s) \right) ds
	\end{equation}
	and
	\begin{equation} \label{maintheorem.pfeq3}
		\|{\bf K}_k(t, s)\|_{{\mathcal S}({\mathbb Z}_k^d, {\mathbb Z}_k^d)}\le  \exp \left({ \|{\bf A}_{k, k}\|_{\infty,{\mathcal S}({\mathbb Z}_k^d, {\mathbb Z}_k^d)}  (t-s)}\right)
		\end{equation}
		for $ t\ge s\ge t_0$. 	By \eqref{Schur.eq1}, 
	\eqref{maintheorem.pfeq2}, \eqref{maintheorem.pfeq3}
	and
	Lemmas \ref{maintheorem.lem1} and \ref{maintheorem.lem2}, we obtain
	\begin{eqnarray}\label{maintheorem.pfeq4}
		\|{\pmb \eta}_{k, N}(t)\|_{\ell^\infty({\mathbb Z}^d_k)} & \hskip-0.01in \le &\hskip-0.01in   \int_{t_0}^t  e^{ \|{\bf A}_{k, k}\|_{\infty, {\mathcal S}({\mathbb Z}_k^d, {\mathbb Z}_k^d)} (t-s)} \left( \sum_{l=k+1}^N \|{\bf A}_{k, l}(t)\|_{\infty, {\mathcal S}({\mathbb Z}_k^d, {\mathbb Z}_l^d)} \right. \nonumber\\
& \hskip-0.01in & \hskip-0.01in
\left. \times ~  \|\pmb \eta_{l, N}(s)\|_{\ell^\infty({\mathbb Z}^d_l)}
		  + \sum_{l=N+1}^\infty  \| {\bf A}_{k,l}\|_{\infty, {\mathcal S}({\mathbb Z}_k^d, {\mathbb Z}_l^d)} \|{\bf z}_l(s)\|_{\ell^\infty({\mathbb Z}_l^d)} \right) ds\nonumber \\
		& \hskip-0.01in \le & \hskip-0.01in  D_0 k  \sum_{l=k+1}^N  R^{k-l-1}
		\int_{t_0}^t  e^{D_0 k (t-s)/R} \|\pmb \eta_{l, N}(s)\|_{\ell^\infty({\mathbb Z}^d_l)} ds\nonumber\\
		& \hskip-0.01in & \hskip-0.01in   +  D_0k \int_{t_0}^t   e^{D_0 k (t-s)/R}  \left(\sum_{l=N+1}^\infty    R^{k-l-1} \|{\bf x}(s) \|_\infty^{l}\right)  ds\nonumber\\
		& \hskip-0.01in \le & \hskip-0.01in  \frac{D_0 k}{R}  \sum_{l=k+1}^N  R^{k-l}
		\int_{t_0}^t  e^{D_0 k (t-s)/R} \|\pmb \eta_{l, N}(s)\|_\infty ds\nonumber\\
		&  \hskip-0.01in  & \hskip-0.01in  + \frac{D_0 k}{R}   \int_{t_0}^t  e^{D_0 k (t-s)/R}  (R-M_0)^{-1} R^{k-N}M_0^{N+1} ds  
	\end{eqnarray}
	for all $t_0\le t\le t_0+T^*$ and  $1\le k\le N$.
	Set
	\begin{equation}\label{maintheorem.pfeq4+} v_k(t) = (R-M_0) R^{N-k} M_0^{-N-1} \|\pmb \eta_{k,N}(t)\|_\infty,\ \ t_0\le t\le t_0+T^*,  \end{equation}
	for $1\le k\le N$.
	Multiplying $(R-M_0) R^{N-k} M_0 ^{-N-1}$ at both sizes of \eqref{maintheorem.pfeq4} yields
	\begin{equation}\label{maintheorem.pfeq5}
		v_k(t)\le   \frac{D_0 k}{R} \int_{t_0}^t e^{D_0 k (t-s)/R} \left(\sum_{l=k+1}^N   v_l(s) +1 \right) ds,\ \  1\le k\le N.
	\end{equation}
	Therefore
	\begin{equation} \label{maintheorem.pfeq6}
		v_2(t)+\cdots+v_N(t)+1\le
		\frac{N^{N-2}}{(N-2)!}  e^{D_0N(t-t_0)/R}, \ t_0\le t\le t_0+T^*,
	\end{equation}
	by  Lemma \ref{maintheorem.lem4}.
	Applying \eqref{maintheorem.pfeq5} with $k=1$ and the estimate in
	\eqref{maintheorem.pfeq6}, we obtain  that
	\begin{equation}
		v_1(t)  \le 
		\frac{ D_0}{R}  \frac{N^{N-2}}{(N-2)!}   \int_{t_0}^t e^{D_0(t-s)/R} e^{D_0N (s-t_0)/R} ds  
		\le \frac{N^{N-2}}{(N-1)!} e^{D_0 N(t-t_0)/R}
	\end{equation}
	hold for all $t\in [t_0, t_0+T^*]$.
	Therefore
	\begin{eqnarray}
		\|{\bf y}_{1, N}(t)-{\bf x}(t)\|_\infty  &\hskip-0.01in  = & \hskip-0.01in \|{\pmb \eta}_{1, N}(t)\|_\infty = \frac{RM_0}{R-M_0} \Big(\frac{M_0}{R}\Big)^N v_1(t)\nonumber\\
		&   \hskip-0.01in  \le  &   \hskip-0.01in  \frac{RM_0}{R-M_0}
		\frac{N^{N-1}}{N!}   \left(\frac{M_0 e^{D_0 (t-t_0)/R}}{R}\right)^N \nonumber\\
		&   \hskip-0.01in  \le  &   \hskip-0.01in  \frac{RM_0}{\sqrt{2\pi}(R-M_0)}
		N^{-3/2}   \left(\frac{e M_0}{R} e^{D_0 (t-t_0)/R}\right)^N ,
		\ \  t_0\le t\le t_0+T^*,
	\end{eqnarray}
	where the last inequality follows from  the Stiriling formula,
	$N!\ge  \sqrt{2\pi} N^{N+1/2} e^{-N}$ for $N\ge 1$.
	This proves the  exponential convergence  
	in  \eqref{maintheorem.eq4}.
\end{proof}

\subsection{Proof of Theorem  \ref{maintheorem3}}
\label{maintheorem2.pfsection}

    To prove Theorem \ref{maintheorem3},
 we need the
following Schur norm estimate for the block matrices
${\bf A}_{k,l}$ for all $k, l\ge 1$ in
\eqref{carleman.eq2}.

\begin{lemma}\label{maintheorem3.lem1}
 If coefficients ${\bf f}_{\pmb\alpha}(t)$ for all ${\pmb \alpha}\in {\mathbb Z}_+^d$ in
the Maclaurin series
 of the vector-valued analytic function ${\bf f}(t, {\bf x})$ in the nonlinear dynamical system  \eqref{dynamicsystem}
 satisfy   Assumptions \ref{assump-1} and \ref{assump-2} and
  \eqref{assumption-3}  for some
 $\nu_0\ge 0$, then
 \begin{equation}  \label{maintheorem3.lem1.eq1}
		{\bf A}_{k,l}(t)=0\ \ {\rm if}\ \  1\le l<k-1,
	\end{equation}
	\begin{equation} \label{maintheorem3.lem1.eq2}
		\sup_{t\ge t_0} \|{\bf A}_{k,l}(t)\|_{{\mathcal S}( {\mathbb Z}^d_k, {\mathbb Z}^d_l)}\le  D_0 k R^{k-l-1} \ \ {\rm if}\  \   1\le k\le l,
	\end{equation}
and
	\begin{equation} \label{maintheorem3.lem1.eq3}
		\sup_{t\ge t_0} \|{\bf A}_{k+1,k}(t)\|_{{\mathcal S}( {\mathbb Z}^d_{k+1}, {\mathbb Z}^d_{k})}\le   \nu_0  k  \ \ {\rm if}\  \   k\ge 1.
	\end{equation}
\end{lemma}

We omit the detailed argument as it is similar to  the one used in the proof of Lemma \ref{maintheorem.lem1}.
 To prove Theorem \ref{maintheorem3}, we next show that the solution
${\bf x}(t), t\ge t_0$,   of the nonlinear dynamical system  \eqref{dynamicsystem} is bounded.

\begin{lemma}\label{maintheorem3.lem2}
Suppose that the solution ${\bf x}(t)$, the initial ${\bf x}_0$, the analytic function ${\bf f}(t,{\bf x})$, and the constant $\epsilon_1$  satisfy assumptions in  Theorem \ref{maintheorem3}. 
	Then,
	\begin{equation} \label{maintheorem3.lem2.eq1}
		\|{\bf x}(t)\|_2\le \max (\|{\bf x}_0\|_2,  R\epsilon_1), \ \ t\ge t_0.
	\end{equation}
\end{lemma}

\begin{proof} By \eqref{nohomogenous.def},
\eqref{assumption-3}, and  Assumptions \ref{assump-1} and \ref{assump-2}, we have
	\begin{eqnarray}\label{maintheorem3.lem2.pfeq1}
		\frac{1}{2} \frac{d \|{\bf x}(t)\|_2^2}{dt} 
& \hskip-0.01in = &  \hskip-0.01in \sum_{j=1}^d f_{j, {\bf 0}}(t) x_j(t)+  \sum_{j=1}^d \lambda_{{\bf e}_j} (x_j(t))^2
		+\sum_{j=1}^d \sum_{|\pmb \alpha|\ge 2} x_j(t) f_{j, \pmb \alpha}(t) {\bf x}_{\pmb \alpha}(t)\nonumber\\
		&  \hskip-0.01in \le &   \hskip-0.01in \nu_0\|{\bf x}(t)\|_\infty  -\mu_0 \|{\bf x}(t)\|_2^2+ \|{\bf x}(t)\|_\infty \sum_{k=2}^\infty \left(\sum_{j=1}^d\sum_{\pmb \alpha\in {\mathbb Z}_k^d} |f_{j, \pmb \alpha}(t)| \right) \|{\bf x}(t)\|_\infty^k
		\nonumber\\
		& \hskip-0.01in \le &  \hskip-0.01in
h(\|{\bf x}(t)\|_2) \|{\bf x}(t)\|_2^2,
	\end{eqnarray}
 where
 $$h(u)=\frac{\nu_0}{u} -\mu_0  + \frac{D_0 u}{R (R-u)}=
 \frac{\nu_0}{u}+\frac{D_0}{R-u}-\frac{D_0+ R\mu_0}{R},\  0<u<R.
$$
Solving the quadratic equation
\begin{equation} \label{maintheorem3.lem2.pfeq2} (1+\eta_1) s^2-(\eta_0+\eta_1) s+ \eta_0=0
\end{equation}
gives $s=\epsilon_0, \epsilon_1$. Hence
 $h(u)>0$ for $u\in (0, R \epsilon_1)$ and $h(u)<0$ for $u\in (R\epsilon_1, R\epsilon_0)$.
	This together with \eqref{maintheorem3.lem2.pfeq1} implies that
	$\frac{d \|{\bf x}(t)\|_2^2}{dt}\le 0$ at a small neighborhood of  $t_0$ when $\|{\bf x}_0\|_2\in (R\epsilon_1, R \epsilon_0)$.
Hence  by \eqref{maintheorem3.eq2},
	$\|{\bf x}(t)\|_2\le \max(\|{\bf x}_0\|_2, R\epsilon_1)$ for all $t_0\le t\le t_0+\delta$ for some $\delta>0$. Using the above argument repeatedly
proves  \eqref{maintheorem3.lem2.eq1}.
%
\end{proof}

For the case that $\nu_0=0$, we have $\epsilon_1=0$. Then applying the argument used in
the proof of Lemma \ref{maintheorem3.lem2}, 
	we obtain
	$$\frac{d \|{\bf x}(t)\|_2^2}{dt}\le -2 \left(\mu_0- \frac{D_0 \|{\bf x}_0\|_2}{R (R-\|{\bf x}_0\|_2)}\right) \|{\bf x}(t)\|_2^2.$$
	Therefore,
	$$\|{\bf x}(t)\|_2\le \|{\bf x}_0\|^2 \exp\left(-2 \left(\mu_0- \frac{D_0 \|{\bf x}_0\|_2}{R (R-\|{\bf x}_0\|_2)}\right)  (t-t_0)\right), \ \ t\ge t_0.$$
	Taking square roots at the above estimate proves the convergence
of the solution ${\bf x}$ of the nonlinear system \eqref{dynamicsystem}   to the equilibrium ${\bf 0}$ exponentially.

\begin{corollary}\label{maintheorem3.cor1}
Suppose that  ${\bf x}(t)$ is  a  solution of  the nonlinear  system \eqref{dynamicsystem} with respect to an initial condition ${\bf x}_0$ that satisfies \eqref{maintheorem2.eq1}. If the  analytic function ${\bf f}(t,{\bf x})$ in   \eqref{dynamicsystem}  meets
	Assumptions \ref{assump-1} and \ref{assump-2}, then
	\begin{equation} \label{maintheorem2.cor1.eq1}
		\|{\bf x}(t)\|_2\le \|{\bf x}_0\|_2 \exp\left(- \left(\mu_0- \frac{D_0 \|{\bf x}_0\|_2}{R (R-\|{\bf x}_0\|_2)}\right)  (t-t_0)\right), \ \ t\ge t_0.
	\end{equation}
\end{corollary}

To prove Theorem \ref{maintheorem3}, we need a technical lemma similar to
Lemma \ref{maintheorem.lem4}.

\begin{lemma}\label{maintheorem3.lem3}
If  $w_1, \ldots, w_N$ are continuous functions  satisfying
	\begin{equation}
		\label{maintheorem3.lem3.eq1}
		0\le w_k(t)\le  \frac{D_0}{R} k \int_{t_0}^t e^{-\mu_0 k (t-s)} \left(\sum_{l=k+1}^{N+1} w_l(s)+\frac{\nu_0}{D_0} w_{k-1}(t)\right) ds\end{equation}
and $w_{N+1}(t)=1$, $w_0(t)=0$,  then
	\begin{equation} \label{maintheorem3.lem3.eq2}
\sum_{l=k}^{N+1} w_l(t)\le \epsilon_0^{k-N-1}
\end{equation}
for all $2\le k\le N$ and $t\ge t_0$,     where $ \nu_0,  \epsilon_0, D_0, R, \mu_0$ are constants in \eqref{assumption-3},  \eqref{maintheorem3.eq3+}
and    Assumptions \ref{assump-1} and \ref{assump-2}.
\end{lemma}

\begin{proof} Take arbitrary $T\ge t_0$ and set $W_{k, T}=\sup_{t_0\le t\le T} w_{k}(t)<\infty$ for $0\le k\le N+1$.
Then it suffices to prove that
\begin{equation}  \label{maintheorem3.lem3.pfeq0}
\sum_{l=k}^{N+1} W_{l, T}\le \epsilon_0^{k-N-1}
\end{equation}
for all $2\le k\le N$.
 By \eqref{maintheorem3.lem3.eq1}, we have
\begin{equation} \label{maintheorem3.lem3.pfeq1}
0\le  W_{k, T}\le  \frac{1}{\eta_1}\sum_{l=k+1}^{N+1} W_{l, T}+ \frac{\eta_0}{\eta_1} W_{k-1, T}, \ 1\le k\le N.
\end{equation}

First,  we prove that
\begin{equation} \label{maintheorem3.lem3.pfeq2}
 W_{k, T}\le   \frac{1-\epsilon_0}{\epsilon_0}
 \sum_{l=k+1}^{N+1} W_{l, T}
\end{equation}
by induction on $1\le k\le N$. 
The conclusion \eqref{maintheorem3.lem3.pfeq2} with $k=1$ follows directly from \eqref{maintheorem3.lem3.pfeq1} with $k=1$, the assumption $W_{0, T}=0$
and the observation that
$(1-\epsilon_0)/{\epsilon_0}\ge {1}/{\eta_1}$
by \eqref{maintheorem3.eq3+}.
Inductively we assume that the conclusion \eqref{maintheorem3.lem3.pfeq2} holds for $1\le k\le N-1$. Then
\begin{equation}  \label{maintheorem3.lem3.pfeq3}
W_{k+1, T}  \le    \frac{1}{\eta_1} \sum_{l=k+2}^{N+1} W_{l, T}+\frac{\eta_0}{\eta_1} W_{k, T}\le \frac{\epsilon_0+\eta_0(1-\epsilon_0)}{\eta_1\epsilon_0} \sum_{l=k+2}^{N+1} W_{l, T} + \frac{\eta_0
(1-\epsilon_0) }{\eta_1 \epsilon_0 }W_{k+1, T}, 
\end{equation}
where the first estimate is obtained from \eqref{maintheorem3.lem3.pfeq1} with $k$ replaced by $k+1$ and the second inequality follows from the inductive hypothesis.
Observe that
$$\eta_1\epsilon_0-\eta_0(1-\epsilon_0)=\frac{\epsilon_0}{1-\epsilon_0} \big(\epsilon_0+\eta_0(1-\epsilon_0)\big)>0
$$
by \eqref{maintheorem3.lem2.pfeq2}.
This together with
\eqref{maintheorem3.lem3.pfeq3}
proves
\begin{equation*}
W_{k+1, T}  \le \frac{\epsilon_0+\eta_0(1-\epsilon_0)}{\eta_1\epsilon_0-\eta_0(1-\epsilon_0)} \sum_{l=k+2}^{N+1} W_{l, T}=
\frac{1-\epsilon_0}{\epsilon_0} \sum_{l=k+2}^{N+1} W_{l, T}.
\end{equation*}
Hence,
 the inductive proof of the statement \eqref{maintheorem3.lem3.pfeq2} can proceed.

 Adding
$\sum_{l=k+1}^{N+1} W_{l, T}$ at both sides of the estimate
\eqref{maintheorem3.lem3.pfeq2}
yields
\begin{equation} \label{maintheorem3.lem3.pfeq4}
 \sum_{l=k}^{N+1} W_{l, T}
 \le   \frac{1}{\epsilon_0}
 \sum_{l=k+1}^{N+1} W_{l, T}.
\end{equation}
Recalling that $W_{N+1, T}=1$ and applying
\eqref{maintheorem3.lem3.pfeq4} repeatedly proves
\eqref{maintheorem3.lem3.pfeq0}.
\end{proof}

We finish this subsection with the proof of Theorem \ref{maintheorem3}.

\begin{proof}  Set $\pmb \eta_{k,N}(t)={\bf y}_{k,N}(t)-{\bf z}_k(t), 1\le k\le N$,
and
 $\pmb \eta_{0,N}(t)=0$.
 Following the argument used in Theorem \ref{maintheorem}, we can show that
	\begin{equation}\label{maintheorem3.pfeq1-}
		\dot{ \pmb \eta}_{k,N}(t)=  \sum_{l=k-1}^N {\bf A}_{k, l}(t) \pmb \eta_{l, N}(t) -\sum_{l= N+1}^\infty {\bf A}_{k,l}(t) {\bf z}_l(t) \end{equation}
	with  zero initial ${\pmb \eta}_{k,N}(t_0)={\bf 0}, 1\le k\le N$.
	By Assumption \ref{assump-2}, one may verify that
	${\bf A}_{k, k}(t)$ is time-independent diagonal matrix  ${\bf A}_k$ with ${\pmb \alpha}$-th diagonal entries  taking value
	$\sum_{j=1}^d \alpha_j \lambda_{{\bf e}_j}\le - k\mu_0  $, where ${\pmb \alpha}=[\alpha_1,\ldots, \alpha_d]$.
	This together with  \eqref{maintheorem3.pfeq1-} implies that
	\begin{equation} \label{maintheorem3.pfeq1}
		\|\exp\big((t-s){\bf A}_k \big)\|_{{\mathcal S}({\mathbb Z}_k^d, {\mathbb Z}_k^d)}\le e^{-\mu_0 k  (t-s)},\  \  t\ge s\ge t_0,\end{equation}
	and
	\begin{eqnarray} \label{maintheorem3.pfeq2}
		{\pmb \eta}_{k, N}(t) & = & \int_{t_0}^t  \exp((t-s){\bf A}_k \big)  \Big( {\bf A}_{k, k-1}(s) {\pmb \eta}_{k-1, N}(s)\nonumber \\
& & \quad + \sum_{l=k+1}^N {\bf A}_{k, l}(s) {\pmb \eta}_{l, N}(s)  -\sum_{l=N+1}^\infty {\bf A}_{k,l}(s) {\bf z}_l(s) \Big) ds.
	\end{eqnarray}
	By \eqref{Schur.eq1},
	\eqref{maintheorem3.pfeq1}, \eqref{maintheorem3.pfeq2}
	and
	Lemmas \ref{maintheorem3.lem1} and  \ref{maintheorem3.lem2}, we obtain
	\begin{eqnarray}\label{maintheorem3.pfeq3}
		\|{\pmb \eta}_{k, N}(t)\|_{\infty} 
		& \hskip-0.01in \le & \hskip-0.01in \nu_0 	k	\int_{t_0}^t  e^{-\mu_0 k (t-s)} \|\pmb \eta_{k-1, N}(s)\|_\infty ds\nonumber\\
& \hskip-0.01in  & \hskip-0.01in +  \frac{D_0 k}{R}  \sum_{l=k+1}^N  R^{k-l}
		\int_{t_0}^t  e^{-\mu_0 k (t-s)} \|\pmb \eta_{l, N}(s)\|_\infty ds
 \nonumber\\
		&  \hskip-0.01in  & \hskip-0.01in  + \frac{D_0 k \left(\max\left\{ \|{\bf x}_0\|_2, R\epsilon_1\right\} \right)^{N+1}}{R  (R-\max \left\{ \|{\bf x}_0\|_2, R\epsilon_1 \right\} )}   \int_{t_0}^t  e^{-\mu_0 k (t-s)}
R^{k-N} ds  
	\end{eqnarray}
	for all $t\ge t_0$ and  $1\le k\le N$.
	Therefore,
	\begin{equation}  \label{maintheorem3.pfeq4}
		w_k(t)\le   \frac{D_0 k}{R} \int_{t_0}^t e^{-\mu_0 k (t-s)} \left(\sum_{l=k+1}^N   w_l(s) +1+ \frac{\nu_0}{D_0} w_{l-1}(s) \right) ds,
	\end{equation}
	where  $w_0(t)=0$ and
	\begin{equation}\label{maintheorem3.pfeq5} w_k(t) =\frac{
 (R-\max \left\{ \|{\bf x}_0\|_2, R\epsilon_1\right\}) R^{N}}{ \left(\max \left\{ \|{\bf x}_0\|_2, R\epsilon_1 \right\} \right)^{N+1}}
 R^{-k}\|\pmb \eta_{k,N}(t)\|_\infty,\   1\le k\le N.  \end{equation}
	Applying Lemma \ref{maintheorem3.lem3} with $k=2$ and
	\eqref{maintheorem3.pfeq4} with $k=1$, we have
	\begin{equation}\label{maintheorem3.pfeq6}
		w_1(t)\le   \frac{D_0}{R} \epsilon_0^{-N+1}
		\int_{t_0}^t e^{-\mu_0 (t-s)}  ds\le \eta_1^{-1} \epsilon_0^{-N+1}.
	\end{equation}
	Therefore
	\begin{eqnarray*}
		\|{\bf y}_{1, N}(t)-{\bf x}(t)\|_\infty  &\hskip-0.01in  = & \hskip-0.01in 
		\frac
{ \left(\max \left\{ \|{\bf x}_0\|_2, R\epsilon_1\right\}\right)^{N+1}}
 {(R-\max \left\{ \|{\bf x}_0\|_2, R\epsilon_1 \right\} ) R^{N-1}}
 w_1(t)\nonumber\\
		&   \hskip-0.01in  \le  &  \hskip-0.01in
\frac{\epsilon_0\max \left\{ \|{\bf x}_0\|_2, R\epsilon_1 \right\}}{ \eta_1 (1-\epsilon_0)}
		\left(\frac{\max \left\{ \|{\bf x}_0\|_2, R\epsilon_1 \right\}}{R\epsilon_0}\right)^N
	\end{eqnarray*}
for all $t\ge t_0$, where the last inequality follow from \eqref{maintheorem3.eq2} and \eqref{maintheorem3.pfeq6}.
\end{proof}

\section{Conclusion and discussions}
\label{conclusion.section}
Several explicit error bounds about convergence of the finite-section approximation of the Carleman linearization of a class of nonlinear systems are presented, where we quantify the time interval over which the convergence happens. When the origin is an asymptotically stable equilibrium of the nonlinear system, it is shown that the convergence holds over the entire time horizon. Furthermore, we show that the convergence over the entire time horizon hold for a nonlinear systems if its drift term, i.e., the zeroth order term in its Maclaurin series, satisfies certain boundedness property.

\bibliographystyle{siamplain}
\bibliography{references}
\end{document}

%% file: ex_shared.tex
\usepackage{lipsum}
\usepackage{amsfonts}
\usepackage{graphicx}
\usepackage{epstopdf}
\usepackage{algorithmic}
\ifpdf
  \DeclareGraphicsExtensions{.eps,.pdf,.png,.jpg}
\else
  \DeclareGraphicsExtensions{.eps}
\fi

\newsiamremark{remark}{Remark}
\newsiamremark{hypothesis}{Hypothesis}
\crefname{hypothesis}{Hypothesis}{Hypotheses}
\newsiamthm{claim}{Claim}

\headers{Carleman Linearization of Nonlinear Systems}{A. Amini, C. Zheng, Q. Sun and N. Motee}

\title{Carleman Linearization of Nonlinear  Systems and Its Finite-Section Approximations
\thanks{ Some preliminary versions of this work were announced in \cite{amini2019carleman, amini2021error}. The authors assert that the content of this manuscript  significantly differs from its conference versions as this work contains several new  and improved results as well as several new case studies with respect to its conference versions.
\funding{
This work was supported in parts by the AFOSR FA9550-19-1-0004, ONR N00014-19-1-2478 and NSF DMS-1816313.} }}

\author{Arash Amini\thanks{Mechanical Engineering and Mechanics, Lehigh University, Bethlehem, PA
  (\email{a.amini@lehigh.edu}, \email{motee@lehigh.edu}).}
\and Cong Zheng\thanks{Department of Mathematics, University of Central Florida, Orlando, FL
  (\email{acongz@knights.ucf.edu}, \email{qiyu.sun@ucf.edu}).}
\and Qiyu Sun\footnotemark[3] \and Nader Motee\footnotemark[2]}

\usepackage{amsopn}

\usepackage{mathrsfs}
\usepackage{amssymb}
\usepackage{amsfonts}
\usepackage{amsbsy}
\usepackage{latexsym}
\usepackage{multirow}
\usepackage{xcolor}
\usepackage{epsfig,psfrag,color}
\usepackage{graphicx}
\usepackage{caption}
\usepackage{subcaption}

\renewcommand{\a}{\alpha}

\renewcommand{\theequation}{\thesection.\arabic{equation}}
\renewcommand{\thefigure}{\thesection.\arabic{figure}}

\newtheorem {theo} {\bf Theorem} [section]

\newtheorem{assum}[theo]{\bf Assumption}

\newcommand{\R}{\mathbb R\,}

\numberwithin{equation}{section}

\newcommand{\x}{{\bf x}}
\newcommand{\xa}{{\bf x}_{\pmb \alpha}}